\documentclass[11pt,a4paper]{amsart}
\usepackage{a4,amssymb,latexsym}
\usepackage[pdfpagemode=None,colorlinks=true,linkcolor=blue,citecolor=blue]{hyperref}

\title[Spectral data for HSL tori]{The spectral data for Hamiltonian stationary 
Lagrangian tori in $\R^4$.}
\author{Ian McIntosh}
\address{Department of Mathematics\\ University of York \\ York Y010 5DD, U.K.}
\email{im7@york.ac.uk}
\author{Pascal Romon}
\address{Laboratoire Analyse et Math\' ematiques Appliqu\' ees\\
Universit\' e Paris-Est Marne-la-Vall\' ee\\
77454 Champs-sur-Marne\\ France}
\email{pascal.romon@univ-mlv.fr}
\date{2nd Dec, 2009}

\newcommand{\N}{\mathbb{N}}
\newcommand{\Z}{\mathbb{Z}}

\newcommand{\C}{\mathbb{C}}
\newcommand{\Ct}{\mathbb{C}^\times}
\newcommand{\R}{\mathbb{R}}
\renewcommand{\H}{\mathbb{H}}

\renewcommand{\P}{\mathbb{P}}
\newcommand{\CP}{\mathbb{CP}}
\newcommand{\HP}{\mathbb{HP}}

\newcommand{\caA}{\mathcal{A}}
\newcommand{\caB}{\mathcal{B}}

\newcommand{\caE}{\mathcal{E}}

\newcommand{\caH}{\mathcal{H}}
\newcommand{\caI}{\mathcal{I}}
\newcommand{\caJ}{\mathcal{J}}
\newcommand{\caK}{\mathcal{K}}
\newcommand{\caL}{\mathcal{L}}
\newcommand{\caM}{\mathcal{M}}
\newcommand{\caO}{\mathcal{O}}

\newcommand{\caS}{\mathcal{S}}
\newcommand{\caT}{\mathcal{T}}
\newcommand{\caU}{\mathcal{U}}

\newcommand{\caX}{\mathcal{X}}
\newcommand{\caZ}{\mathcal{Z}}

\newcommand{\fS}{\mathfrak{S}}
\newcommand{\so}{\mathfrak{so}}
\newcommand{\gl}{\mathfrak{gl}}
\newcommand{\fg}{\mathfrak{g}}
\newcommand{\fp}{\mathfrak{p}}

\newcommand{\ft}{\mathfrak{t}}
\newcommand{\fz}{\mathfrak{z}}
\newcommand{\Lgt}{\Lambda^\tau\fg^\C}
\newcommand{\Lg}{\Lambda^\tau\fg}

\newcommand{\Lpm}{\Lambda^\mu\fp}
\newcommand{\Lp}{\Lambda\fp}

\newcommand{\Lag}{\mathrm{Lag}}
\newcommand{\Ham}{\mathrm{Ham}}
\newcommand{\End}{\mathrm{End}}
\newcommand{\Pic}{\mathrm{Pic}}
\newcommand{\Div}{\mathrm{Div}}
\newcommand{\Jac}{\mathrm{Jac}}

\newcommand{\Spec}{\mathrm{Spec}}
\newcommand{\Spin}{\mathit{Spin}}

\renewcommand{\Re}{\mathrm{Re}}

\renewcommand{\Im}{\mathrm{Im}}

\newcommand{\Kah}{K\" ahler\ }
\newcommand{\pk}{polynomial Killing\ }
\newcommand{\HSL}{HSL\ }
\newcommand{\Ad}{\mathrm{Ad}}
\newcommand{\bt}{\mathbf{t}}
\newcommand{\ba}{\mathbf{a}}
\newcommand{\e}{\varepsilon}

\newtheorem{thm}{Theorem}[section]
\newtheorem{prop}[thm]{Proposition}
\newtheorem{lem}[thm]{Lemma}
\newtheorem{cor}[thm]{Corollary}

\theoremstyle{remark}

\newtheorem{rem}[thm]{Remark}

\begin{document}

\begin{abstract}
This article determines the spectral data, in the integrable systems sense, 
for all weakly conformally immersed Hamiltonian stationary Lagrangian in $\R^4$.
This enables us to describe their moduli space and the locus of branch points of
such an immersion. This is also an informative example in integrable systems 
geometry, since the group of ambient isometries acts non-trivially on the spectral data 
and the relevant energy
functional (the area) need not be constant under deformations by higher flows.
\end{abstract}

\maketitle

\section{Introduction.}

A smooth immersion of a surface $f:M\to\R^4$ is Hamiltonian stationary Lagrangian\footnote{The original
terminology for these was H-minimal or Hamiltonian minimal.} (HSL) if $f(M)$ is a
Lagrangian submanifold whose area is
stationary for all variations by (compactly supported) Hamiltonian vector fields.
This is a natural variational problem for Lagrangian submanifolds and occurs in the study of
volume minimisers in families of Lagrangian submanifolds \cite{SchW,SchW1}. The Euler-Lagrange
equations were derived by Oh \cite{Oh2} and can be phrased in terms of the 
the mean curvature 1-form $\sigma_H=f^*(H\rfloor\omega)$, where $H$ denotes the mean curvature vector and 
$\omega$ is the standard K\" ahler form on $\R^4\simeq\C^2$. 
For any Lagrangian submanifold of $\R^{2n}$ it can be shown that $d\sigma_H=0$, and $f$ is 
Hamiltonian stationary
when $d*\sigma_H=0$, i.e., $\sigma_H$ is harmonic. 
Consequently, for a compact HSL surface $\sigma_H$ represents a cohomology class. It turns out
that, after scaling, this is the Maslov class $\mu\in H^1(M,\Z)$, which is an important Hamiltonian
isotopy invariant of $f(M)$.
Oh conjectured in \cite{Oh2} that the Clifford torus in
$S^3\subset\R^4$ minimises area in its Hamiltonian isotopy class and this proved to be a challenging
question, which remains unanswered (although see \cite{Ilm, Anc,Anc2} for a partial solution).

Rather surprisingly, all weakly conformal HSL immersions of a torus $\C/\Gamma$ into $\R^4$ can be
explicitly described in terms of the Fourier components of the immersion: this was discovered by
H\' elein \& Romon \cite{HelR}, who showed that the equations can be reduced to a system of 
linear equations. 
They arrived at this through an investigation of the HSL equations as the Maurer-Cartan equations for a
loop of flat connexions, where the connexion form takes values in the Lie algebra of the group $G$ of
symplectic isometries of $\R^4$. This is an ``integrable systems'' approach, and although this proved
to be somewhat superfluous to writing down the Fourier component solution, H\' elein \& Romon
showed that the problem does admit
the structure normally associated with an integrable system. The loop algebra valued Maurer-Cartan form
is quadratic in the loop parameter,
\begin{equation}
\alpha_\zeta = \zeta^{-2}\alpha_{-2} + \zeta^{-1}\alpha_{-1} + \alpha_0
+ \zeta\alpha_1 +\zeta^2\alpha_2,
\end{equation}
and satisfies a Lax equation, $d\alpha_\zeta = [\xi_\zeta,\alpha_\zeta]$ where $\xi_\zeta$ is a Laurent
polynomial in $\zeta$ (that is to say, HSL tori possess a ``\pk field''). However, the Lax
equations have a significant difference from those which are most often seen in the integrable systems
study of, say, harmonic maps or minimal surfaces. The group $G$ is not semisimple (or even reductive):
it is the semi-direct product $U(2)\ltimes\C^2$ and its complexification is most naturally realised
inside a parabolic subgroup of $GL_5(\C)$. Together with the fundamental linearity of this particular
set of equations, this raises the question of what the solution to these Lax equations look like, and
particularly whether there is any effective ``spectral data'' which reproduces the Fourier component
solutions. Here we are referring to the expectation of a correspondence
between a HSL torus and algebro-geometric
data of the type which one typically sees in the solution of Lax equations.

Our aim here is to show that such a correspondence does exist and not only does it reproduce entirely
the solutions from the Fourier component method, it provides extra information which that method does
not yield so easily. In particular, given that the original formulation allows for branch points it
would be nice to know to what extent branched immersions make up the vector space of solutions found in
\cite{HelR}. Our approach leads to the conclusion that unbranched immersions are 
generic, and generically amongst branch immersions the branch points are isolated (see remark
\ref{rem:branchpts}). 

We first show that every weakly conformal HSL immersion
$f:\C/\Gamma\to\R^4$ is determined by a triple of data $(X,\lambda,\caL)$ consisting of a  complete
algebraic curve $X$, a rational function $\lambda$ on $X$ and a line bundle $\caL$ (or possibly a
non-invertible rank $1$ sheaf) over $X$. To achieve this we cannot use 
just the naive characteristic polynomial spectral curve of one \pk field $\xi_\zeta$, since this does
not yield enough information. Instead we use a commutative subalgebra of \pk fields. 
The linear nature of the problem forces $X$ to be a rational
(and in this case, reducible) curve: its singularities correspond precisely to the non-trivial Fourier
modes of $f$. More suprisingly, the parabolic nature of the group $G$ means the spectral data is
\emph{not} invariant under symplectic isometries. This is a consequence of the fact  
not all commutative subalgebras of \pk fields are isomorphic, so that the choice of subalgebra becomes
part of the spectral data. 

When this is all put together the full picture has a simple geometry. 
Fix a conformal type $\Gamma$ for the torus and fix a Maslov
class $\beta_0\in\Gamma^*\simeq H^1(\C/\Gamma,\Z)$. This fixes the spectral curve $(X,\lambda)$.
Let $\caS(\Gamma,\beta_0)$ be the set of triples $(X,\lambda,\caL)$
for weakly conformal HSL immersions $f:\C/\Gamma\to\R^4$ of Maslov class $\beta_0$ satisfying the base
point condition $f(0)=0$ (and a similar condition on the Lagrangian angle function). This 
condition is preserved by dilations, which always preserve the HSL condition, and a subgroup $G_0$ of
symplectic isometries isomorphic to $SU(2)$. Let $N$ denote the maximum number of non-trivial Fourier
modes such a map can possess. We prove:
\begin{thm}
$\caS(\Gamma,\beta_0)\simeq\CP^{N-1}$ and corresponds to the space of based weakly conformal HSL immersions
$f$ modulo dilations and the action of a maximal torus $S^1\subset G_0$. Varying the base point of $f$
induces an action of the real Lie group $\C/\Gamma$ on
$\caS(\Gamma,\beta_0)$ through a map $p\mapsto (X,\lambda,\caL_p)$. 
The map $f$ has a branch point at $p\in\C/\Gamma$ 
precisely when $\caL_p$ lies on the intersection of
two hyperplanes $\Theta_\infty$ and $\Theta_0$ in $\caS(\Gamma,\beta_0)$.
\end{thm}   
The quotient of
$\caS(\Gamma,\beta_0)$ by the action of $G_0$ is a $S^2$-bundle over the moduli space
$\caM(\Gamma,\beta_0)$ of these HSL tori. 
The natural map $\caS(\Gamma,\beta_0)\to\caM(\Gamma,\beta_0)$ can be realised as the fibration
$\CP^{N-1}\to\HP^{N/2-1}$, and the fibres are the $G_0$-orbits of spectral data. 
We provide an explicit expression for the hyperplanes $\Theta_\infty,\Theta_0$ in the homogeneous 
coordinates
on $\caS(\Gamma,\beta_0)$. Each can be interpreted as a translate of the $\theta$-divisor for
the $\theta$-function of a rational curve: this is explained in the appendix. 
It follows that the parameterisation by spectral data turns the question of locating branch points 
into a problem in projective geometry.

As one expects, the Jacobi variety of $X$, or more precisely a real subgroup $J_R$ of it, acts on the
spectral data to produce new HSL tori. The Maslov class alone determines $(X,\lambda)$ completely, so
that $\caS(\Gamma,\beta_0)$ is a union of $J_R$-orbits of different dimensions. The sheaf $\caL$ is
only a line bundle on the largest orbit. The base point translation action of $\C/\Gamma$ factors
through a homomorphism of $\C/\Gamma$ into $J_R$. In integrable systems language the action of $J_R$
generates the ``higher flows''. We show in \S 6 that these do not all correspond to Hamiltonian 
variations.
This too is a surprising departure from the study of minimal tori, where all higher flows preserve
the area. It reflects the fact that the area functional is a non-constant function on the moduli
space $\caS(\Gamma,\beta_0)$.

We are hopeful this study will provide a
useful testing ground for deeper investigations in the theory of spectral data in surface theory.
In particular, the biggest challenge at this time is to understand how the appropriate energy 
functional depends on the spectral genus. For HSL tori in $\R^4$ the spectral data is
elementary enough to see (cf.\ remarks \ref{rem:area} and \ref{rem:genus} below) how the area 
functional depends upon 
the conformal class and the spectral data.  We hope this relationship which may provide 
insight into how this works more generally in integrable surface theory. 

\section{Lagrangian surfaces in $\R^4$.}

Let $\R^4$ be equipped with its Euclidean metric, its standard complex structure $J$, and 
\Kah form $\omega$. We will represent its group of symplectic isometries in the form
\[
G = \{(g,u)\in SO(4)\ltimes\R^4: gJg^{-1} = J\}\simeq U(2)\ltimes\C^2.
\]
Now suppose we have a conformally immersed orientable 
surface $f:M\to\R^4$. If it is
Lagrangian (i.e., $f^*\omega=0$) then its Gauss map $\gamma:M\to\Lag(\R^4)$ takes values in the
Grassmannian $\Lag(\R^4)$ of oriented Lagrangian $2$-planes in $\R^4$.
Let $H:M\to TM^\perp$ be the mean curvature field for $f$, then one knows
that the mean curvature form $\sigma_H=f^*(H\rfloor\omega)$ is closed and $f$ is
Hamiltonian stationary precisely when $\sigma_H$ is also co-closed. Now,
since $\Lag(\R^4)\simeq U(2)/SO(2)$ we can post-compose $\gamma$ with the
well-defined map $\det:U(2)/SO(2)\to S^1$ induced by the determinant on
$U(2)$. Let $s:M\to S^1$ be defined by $s = \det\circ\gamma$, then by definition the
Maslov form $\mu\in\Omega^1_M$ of $f$ is 
\[
\mu = \frac{1}{\pi i}s^{-1}ds.
\]
By a theorem of Morvan \cite{Mor} this is related to the mean curvature form by $\mu =
\frac{2}{\pi}\sigma_H$ and therefore
\begin{equation}\label{eq:s}
s^{-1}ds = 2i\sigma_H.
\end{equation}
Thus $f$ is Hamiltonian stationary if and only if $s$ is a harmonic
map. 

Now we restrict our attention to
the case where $M$ is a torus, represented in the form $M=\C/\Gamma$
where $\Gamma$ is a lattice. We can write $s = \exp(i\beta)$ for a
function $\beta:\C\to\R$, which is called the Lagrangian angle, and $f$ is Hamiltonian stationary when $\beta$
is a harmonic function. 
Note that the action of the centre of $U(2)$ on $f$ is by $f\mapsto e^{\theta J}f$ (for some
$\theta\in[0,2\pi]$), under which the Lagrangian angle changes by $\beta\mapsto\beta+2\theta$. Therefore
we may (and will) assume that $\beta(0)=0$. 
Following \cite{HelR}, and using the inner product $\langle z,w\rangle =
\Re(z\bar w)$ on $\C$, we can write $\beta$ as 
\begin{equation}\label{eq:beta}
\beta(z) = 2\pi\langle \beta_0,z\rangle 
\end{equation}
for a constant $\beta_0\in\Gamma^*\subset\C$, i.e., $\beta(z)\in 2\pi\Z$ for every $z\in\Gamma$.
\begin{rem}
The Maslov class of $f$ is the cohomology class of $[\mu]\in H^1(M,\Z)$. Under the natural
identification $H^1(M,\Z)\simeq\Gamma^*$ we can think of $\beta_0$ as the Maslov class.
One knows that \HSL surfaces are constrained Willmore surfaces \cite{Boh08} (i.e., the Willmore energy
$W(f)=\int_M|H|^2$ is critical for variations through conformal immersions). Because of the relation $H
= -\frac{1}{2}J\nabla\beta$ the Willmore energy for a conformally immersed HSL torus $f:\C/\Gamma\to\R^4$ 
is ``quantized'' by the
Maslov class $\beta_0$:
\begin{equation}\label{eq:Willmore}
W(f) = \pi^2|\beta_0|^2A(\C/\Gamma),
\end{equation}
where $A(\C/\Gamma)$ is the area of the flat torus $\C/\Gamma$ with metric $|dz|^2$.
\end{rem}

\subsection{Twistor lift and frames.}
Let $E$ denote the pullback $f^{-1}T\R^4$ of the tangent
bundle of $\R^4$. Then $E=TM\oplus TM^\perp$ and since $f$ is Lagrangian $TM^\perp = JTM$. Let
$J_M$ denote the intrinsic complex structure carried by $M$. Then since $f$ is conformal it
induces another complex structure $S=J_M\oplus JJ_MJ$ on $E$, with the property that $JS=-SJ$. We
may think of $S$ as a twistor lift \cite{BurK}
of $f$, i.e., $S:M\to Z$ where $Z$ is the twistor 
bundle of complex structures on $T\R^4$. Inside $Z$ lies the $S^1$-subbundle of all complex 
structures which anti-commute with $J$, which is where $S$ takes values. In fact this
$S^1$-bundle is the image of a 4-symmetric space $G/G_0$, where $G_0$ is the fixed point
subgroup of an order 4 outer automorphism $\tau$ of $G$.
To see this, first let $\e_1,\ldots,\e_4$ be the standard oriented orthonormal basis of $\R^4$ 
for which $\e_2=J\e_1$, $\e_4=J\e_3$, and 
define $L\in SO(4)$ to be the complex structure on $\R^4$ characterised
\[
L\e_1=\e_3,\ L\e_2=-\e_4. 
\]
We observe that $LJ=-JL$ and that every complex structure anti-commuting with $J$ is of the form
$gLg^{-1}$ for some $g\in U(2)\subset SO(4)$. Using $L$ we define an order 4 outer automorphism 
\[
\tau:G\to G; \ \tau(g,u) = (-LgL,-Lu).
\]
The fixed point subgroup is $G_0= \{(g,0)\in G: gL=Lg\}\simeq SU(2)$ and the 4-symmetric space
$G/G_0$ is an $S^1$-bundle over $\R^4$. On the other hand the twistor
bundle has description 
\[
Z=\{(gJg^{-1},u)\subset SO(4)\times \R^4:g\in SO(4)\},
\]
which is isomorphic to the homogeneous space $(SO(4)\ltimes\R^4)/U(2)$. Now we have the embedding
\[
G/G_0\to Z;\ (g,u)G_0\mapsto (gLg^{-1},u),
\]
whose image is the bundle of complex structures anti-commuting with $J$. 

Let us now give a
description of $S:M\to Z$ in terms of natural frames, and show that $S$ is essentially the
Lagrangian angle function.
Any conformal Lagrangian torus possesses a natural frame on the universal cover $\C$, 
called the fundamental frame, $\tilde{f}:\C\to G$, given by
$\tilde{f} = (F,f)$, where $F:\C\to SO(4)$ is chosen so that $F\e_j=f_j$, where
\[
f_1 = e^{-\rho}f_x,\ f_2 = e^{-\rho}Jf_x,\ 
f_3 = e^{-\rho}f_y,\ f_4 = e^{-\rho}Jf_y, 
\]
for $|df|^2 = e^{2\rho}|dz|^2$, where $z=x+iy$. Setting  $\epsilon = (\e_1-i\e_3)/2$ we see
\[
df = e^\rho F(\epsilon dz +\bar\epsilon d\bar z).
\]
It follows that $s = \det(F)$.\footnote{Beware here that $\det(F)$ is the
determinant taken in $U(2)$ not in its representation as a subgroup of
$O(4)$.} 

From the definition of $S$ in terms of the complex structures $J_M$ and $J$ it follows that
$Sf_1 = f_3$ and $Sf_2 = -f_4$. Therefore $S = (FLF^{-1},f)$. But, as observed in \cite{HelR},
we can also take one of two spinor frames for $f$:
\[
U_{\pm}= (\pm\exp(J\beta/2),f).
\]
Notice that $\tilde{f} = U_\pm K$ where $K = (\pm\det(F)^{-1/2}F,0)$. Since $K$ takes values in
$G_0$ we have 
\[
FLF^{-1} = e^{J\beta/2}Le^{-J\beta/2} = e^{J\beta}L,
\]
and therefore $S$ is essentially the Lagrangian angle function. 
\begin{rem}\label{rem:S-hol}
This gives us another perspective on the Hamiltonian stationary condition, namely, a map $f:M\to\R^4$
is conformal Lagrangian if and only if it is $S$-holomorphic for some $S:M\to G/G_0$ (i.e., $Sf_z =
if_z$). In that case we necessarily we have $S=e^{\beta J}L$ for some function $\beta$. Then 
$f$ is Hamiltonian stationary if
further $\Delta\beta = 0$. 
\end{rem}

Throughout the remainder of this article we choose to work with the spinor lift $U_+$ since it
has a particularly nice Maurer-Cartan form: we will define 
\[
\alpha = U_+^{-1}dU_+ = (\frac{1}{2}Jd\beta, e^{-J\beta/2}df).
\]
We also assume, without any loss of generality, that $f(0)=0$ and therefore there is a unique spinor
lift $U_+$ determined by $\alpha$ satisfying the initial condition $U_+(0)=(I,0)$.
The subgroup of symplectic isometries preserving the two conditions $f(0)=0$ and $\beta(0)=0$ is $G_0$.

\subsection{Extended Maurer-Cartan form.} 
The automorphism $\tau$ induces an order 4 automorphism (which we
shall also call $\tau$) on $\fg^\C$, the complexification of the Lie algebra
$\fg$ of $G$.  We will represent $\fg^\C$, as vector space, by
\[
\fg^\C = \{(X,x)\in\so_4(\C)\times \C^4: [X,J]=0\}.
\]
The automorphism $\tau$ takes the form
\[
\tau(X,x) = (-LXL,-Lx),\ 
\]
Let $\fg_j\subset\fg^\C$ be the $i^j$-eigenspace for $\tau$, then one computes
\begin{eqnarray*}
\fg_{-1} = \C\epsilon\oplus\C J\bar\epsilon,&  \fg_0 = \{(X,0)\in\fg^\C): [X,L]=0\},\\
 \fg_1 = \bar\fg_{-1},& \fg_2 = \{(rJ,0):r\in\C\}.
\end{eqnarray*}
Notice that $\fg_{-1}$ is the $i$-eigenspace for $L\in\End(\C^4)$.
We notice that the components of $\alpha$ in this decomposition are 
\[
\alpha_{-1} = (0,e^{-J\beta/2}\frac{\partial f}{\partial z}dz),\ 
\alpha_{1} = (0,e^{-J\beta/2}\frac{\partial f}{\partial \bar z}d\bar z),\ 
\alpha_0=(0,0),\ \alpha_2 = (\frac{1}{2}Jd\beta,0).
\]
We define the extended Maurer Cartan form to be the loop of 1-forms
\begin{equation}\label{eq:alpha}
\alpha_\zeta = \zeta^{-2}\alpha_2^\prime + \zeta^{-1}\alpha_{-1} + \alpha_0
+ \zeta\alpha_1 +\zeta^2\alpha_2^{\prime\prime}.
\end{equation}
It is the principal observation of H\' elein \& Romon \cite{HelR} that the Maurer-Cartan equations 
for $\alpha_\zeta$
are satisfied if and only if $f$ is Hamiltonian stationary. As usual, we think of $\alpha_\zeta$ as a
1-form with values in a loop algebra. It possesses two symmetries, namely, a real symmetry and
$\tau$-equivariance:
\begin{equation}\label{eq:real}
\overline{\alpha_{\zeta}} = \alpha_{\bar\zeta^{-1}},\quad
\tau(\alpha_\zeta) = \alpha_{i\zeta},
\end{equation}
where $\overline{(X,x)} = (\bar X,\bar x)$ is simply complex conjugation. 
Therefore we may think of $\alpha_\zeta$ as taking values in
the twisted loop algebra $\Lg$ of $\tau$-equivariant real analytic maps $\xi_\zeta:S^1\to\fg^\C$
possessing the real symmetry. This is a real subalgebra of the complex algebra $\Lgt$ of
$\tau$-equivariant real analytic maps $\xi_\zeta:S^1\to\fg^\C$.

\section{Polynomial Killing fields.}

H\' elein \& Romon \cite{HelR} have shown that every Hamiltonian stationary
Lagrangian torus in $\R^4$ has an adapted \pk field, i.e., a
map $\xi_\zeta:\C/\Gamma\to \Lg$ satisfying 
\begin{enumerate}
\item $d\xi_\zeta = [\xi_\zeta,\alpha_\zeta]$,
\item $\xi_\zeta = \zeta^{-4d-2}\alpha_{-2} + \zeta^{-4d-1}\alpha_{-1} + \ldots.$
\end{enumerate} 
However, there are infinitely many linearly independent adapted \pk fields. 
Following the principle in \cite{McI} we would like to say that, by dropping condition (b) and allowing
$\xi_\zeta$ to take values in $\Lgt$,
\pk fields come
in complex algebras, and that the solution of the Lax equation (and the geometry of the original map) should
be able to be reconstructed from spectral data determined by this algebra. However, 
matrix multiplication does not preserve the loop algebra $\Lgt$. The reason is
that $\fg$ itself is not closed under matrix multiplication.
We rectify this by working in the larger matrix algebra
\[
\fp = \{\left(\begin{smallmatrix} A & a\\ 0 & b\end{smallmatrix}\right)\in \gl_5(\C):
A\in\gl_4(\C), [A,J]=0,a\in\C^4,b\in \C\}.
\]
This contains $\fg^\C$ as the subalgebra for which $A\in\so_4(\C)$ and $b=0$.
Rather than use the block matrix notation, it will be more convenient to use
a hybrid notation,
\begin{equation}\label{eq:hybrid}
\left(\begin{smallmatrix} B & a\\ 0 & b\end{smallmatrix}\right)=(A,a) + bI,
\end{equation}
where $I$ stands for the identity matrix in $\gl_5$ and we have defined
$A = B-bI_4$. In this notation the Lie bracket is
given by
\[
[(X,x)+yI,(A,a)+bI] = ([X,A],Xa-Ax).
\]
Notice that $\tau$ extends to $\fp$ as the Lie algebra automorphism for which
\[
\tau:(X,x)+yI\mapsto (-LXL,-Lx)+yI.
\]
To properly understand the algebra of \pk fields and their spectral data, we will 
work in a different realisation of the twisted loop algebra, in which the 
twisting is partially removed. First we embed $\Lgt$ in
$\Lpm$, where $\mu=\tau^2$. The involution $\mu$ is inner and
therefore $\Lpm$ is an algebra with unit under matrix multiplication. 
To remove the twisting in $\Lpm$ write $\mu=\Ad Q$ 
and let $\kappa_\zeta:S^1\to K$ be a homomorphism 
for which $\kappa_1=I$ and $\kappa_\omega = Q^{-1}$. 
Since $\mu:\fg^\C\to\fg^\C:(A,a) \mapsto (A,-a)$, we have
\[
\Ad\kappa_\zeta\cdot (X_\zeta,x_\zeta) = (X_\zeta,\zeta x_\zeta).
\]
The right hand side is an untwisted loop when considered as a function of $\lambda=\zeta^2$. This
extends naturally to $\Lpm$ to give a matrix algebra automorphism from $\Lpm$ to the algebra 
$\Lp$ of untwisted loops in $\fp$. From now on 
we shall assume this has been applied to all the objects under study:  
the notation will imply that the untwisting has been applied by writing all loops as a function of
$\lambda=\zeta^2$. 

In $\Lp$ the extended Maurer-Cartan form has the shape
\begin{equation}\label{eq:untwisted}
\alpha_\lambda =  (\frac{\pi}{2}(\lambda^{-1}\bar\beta_0 dz+\lambda\beta_0d\bar z)J, 
e^{-J\beta/2}(f_zdz+\lambda f_{\bar z}d\bar z)).
\end{equation}
Since $\alpha_\zeta$ is both $\tau$-equivariant and satisfies the reality condition
\eqref{eq:real}, this form
$\alpha_\lambda$ has the induced symmetries
$\bar\rho^*(\alpha_\lambda) = \alpha_\lambda$, $\tau^*(\alpha_\lambda) = \alpha_\lambda$, where 
$\bar\rho^*,\tau^*$  are the commuting, respectively $\R$-linear and $\C$-linear, algebra involutions of $\Lp$ 
defined by 
\begin{equation}\label{eq:involutions}
\bar\rho^*(\xi_\lambda)  =  \Ad R_\lambda\cdot \overline{\xi_{\bar\lambda^{-1}}},
\quad \tau^*(\xi_\lambda) = \Ad T\cdot\tau(\xi_{-\lambda}),
\end{equation}
where 
\[
R_\lambda = \begin{pmatrix} I_4 & 0\\ 0 & \lambda^{-1}\end{pmatrix},\ 
T = \begin{pmatrix} iI_4 & 0\\ 0 & 1\end{pmatrix}.
\]
For $\xi_\lambda = (X_\lambda,x_\lambda) + y_\lambda I$ these look like
\begin{equation}\label{eq:rhotau^*}
\bar\rho^*(\xi_\lambda) = (\bar X_{\bar\lambda^{-1}},\lambda \bar x_{\bar\lambda^{-1}}) + \bar
y_{\bar\lambda^{-1}}I,\quad
\tau^*(\xi_\lambda) = (-LX_{-\lambda}L,-iLx_{-\lambda}) + y_{-\lambda}I.
\end{equation}
Both $\bar\rho^*$ and $\tau^*$ are loop algebra automorphisms.
As a result of this symmetry of $\alpha_\lambda$, if $\xi_\lambda$ is a \pk field then so are both
$\bar\rho^*(\xi_\lambda)$ and $\tau^*(\xi_\lambda)$. 

From now on we extend the definition of a \pk field to include any solution of 
(a) with values in $\Lp$ .  Our aim is to understand the algebra $\caK$ of these \pk fields
and, more precisely, to describe certain maximal abelian subalgebras.
For a \pk field of the form $(X_\lambda,x_\lambda) + y_\lambda I$ the \pk field equations 
become 
\begin{eqnarray}\label{eq:pkf1}
dX_\lambda & = & 0 \notag \\
dy_\lambda & = & 0 \notag\\
dx_\lambda + \frac{\pi}{2}(\lambda^{-1}\bar\beta_0 dz+\lambda\beta_0d\bar z)Jx_\lambda& 
=& X_\lambda e^{-J\beta/2}(f_zdz+\lambda f_{\bar z}d\bar z).
\end{eqnarray}
Therefore $X_\lambda$ and $y_\lambda$ depend on $\lambda$ alone. Since these equations are linear over
$\C[\lambda,\lambda^{-1}]$ we may assume that $X_\lambda$ is a polynomial of degree $N$, $X_\lambda =
\sum_{j=0}^NX_j\lambda^j$. Expanding $x_\lambda$ similarly we obtain the equations
\begin{eqnarray}
\partial x_j/\partial z + \frac{\pi}{2}\bar\beta_0Jx_{j+1}& =& X_je^{-J\beta/2}f_z,\label{eq:pkf2}\\
\partial x_j/\partial\bar z + \frac{\pi}{2}\beta_0Jx_{j-1}& =& X_{j-1}e^{-J\beta/2}f_{\bar
z}.\label{eq:pkf3}
\end{eqnarray}
It follows that $x_j=0$ for $j\leq 0$ and $j>N$. It is also clear that if $X_\lambda$ is identically
zero then so is $x_\lambda$, i.e., there are no non-trivial \pk fields of the form $(0,x_\lambda)$. 
Now consider \pk fields for which $X_\lambda = q_\lambda R$, where $q_\lambda\in\C[\lambda]$
and $R$ is a constant projection matrix (i.e., $R^2=R$), which we will assume to be orthogonal
projection onto some subspace $V$ of $\C^4$: since $[R,J]=0$ this subspace $V$ is 
$J$-invariant.\footnote{We know such \pk fields must exist since in \cite{HelR} it is shown that 
there are \pk fields 
of the form $(p(\lambda) J,y_\lambda)$ and multiplication on the left by $-RJ$ produces solutions of
\eqref{eq:pkf1} by linearity.}
\begin{lem}\label{lem:minimal}
For each orthogonal projection matrix $R$ commuting with $J$
there exists a unique monic polynomial $p(\lambda)$, of minimal degree $N$ with $p(0)\neq 0$,
such that there is a \pk field of the form $(p(\lambda) R, x_\lambda)$. Further, for any 
\pk field of the form $(q(\lambda) R,y_\lambda)$, with $q(\lambda)\in\C[\lambda,\lambda^{-1}]$,
there exists $r(\lambda)\in\C[\lambda,\lambda^{-1}]$ for which
$q(\lambda)= r(\lambda) p(\lambda)$ and $y_\lambda = r(\lambda) x_\lambda$. Finally, $x_\lambda$ takes
values in the image of $R$.
\end{lem}
\begin{proof}
Suppose $(p(\lambda) R,x_\lambda),(p'(\lambda)R,x_\lambda')$ are both non-trivial \pk fields for which 
$p(\lambda),p'(\lambda)$ are monic of the same degree, minimal for \pk fields of this type.  
Then $((p(\lambda)-p'(\lambda))R,x_\lambda-x_\lambda')$ is again a \pk field and so 
$p(\lambda)=p'(\lambda)$, otherwise the degree of their difference is less than the minimal degree.
It follows that $x_\lambda = x_\lambda'$.

Now suppose $(q(\lambda) R,y_\lambda)$ is a \pk field, then there exists $k\in\N$ for which
$\lambda^k q(\lambda)$ is a polynomial. Clearly we can always find a polynomial $s(\lambda)$ 
which $\lambda^k q(\lambda)-s(\lambda) p(\lambda)$ has degree less than $p(\lambda)$. But the 
\pk field equations are linear over $\C[\lambda,\lambda^{-1}]$, therefore 
must have
\[
\lambda^k q(\lambda)-s(\lambda) p(\lambda)=0,\ \lambda^k y_\lambda-s(\lambda) x_\lambda=0.
\]
Finally, $x_\lambda = Rx_\lambda +R^\perp x_\lambda$, where $R^\perp$ is the complementary orthogonal
projection onto $U^\perp$. Since $[R^\perp,J]=0$, $R^\perp x_\lambda$ satisfies \eqref{eq:pkf1} with
$X_\lambda = 0$, hence it is identically zero. 
\end{proof}
For convenience, we will refer to the \pk field with this unique, minimal degree, monic polynomial 
multiplier of $R$ as the minimal \pk field for the projector $R$.
Note that if $\dim(V)=1$ then $x_\lambda = s_\lambda v$ for some non-zero $v\in V$ and some function
$s_\lambda(z)$ which is polynomial in $\lambda$.

Everything we need to know about the spectral curve is encoded in the minimal
\pk field of the form $(p(\lambda) I_4,x_\lambda)$.
From now on
we will use $p(\lambda),x_\lambda$ exclusively for this \pk field, and denote the degree of $p(\lambda)$
by $N$. To understand $p(\lambda)$ we note first that
$x_\lambda$ is completely determined by $p(\lambda)$, given $u$, using the recursion implicit 
in the equations \eqref{eq:pkf2}: we have 
\begin{equation}\label{eq:recursion}
x_{j+1}  = 
-(\frac{2}{\pi\bar\beta_0}J)^{j+1}\sum_{k=0}^j(-\frac{\pi\bar\beta_0}{2}J)^kp_k\frac{\partial^{j-k}
u}{\partial z^{j-k}},\quad 0\leq j\leq N,
\end{equation}
where for simplicity we set $u := \alpha_{-1}(\partial/\partial z) = e^{-J\beta/2}f_z$
(following the notations in \cite[\S 3.2]{HelR}).
Recall also that the Maurer-Cartan equation for $\alpha_\zeta$ imply in particular
\begin{equation}\label{eq:linear}
u_{\bar z} = \frac{\pi}{2}\bar\beta_0J\bar u,\\
\end{equation}
This recursion has an important consequence for the roots of $p(\lambda)$.
Indeed, considering
\eqref{eq:recursion} for $j=N$, we see that, since $x_{N+1}=0$,
\begin{equation}\label{eq:psum}
\sum_{k=0}^N (-\frac{\pi\bar\beta_0}{2}J)^kp_k\frac{\partial^{N-k} u}{\partial z^{N-k}} = 0.
\end{equation}
Now use the Fourier expansions of $f$ and $u$, 
\[
f=\sum_{\gamma\in\Gamma^*} f_\gamma e_\gamma,\quad u=\sum_{\gamma\in\frac{1}{2}\Gamma^*} u_\gamma
e_\gamma,
\]
where $e_\gamma(z)= \exp(2\pi i \langle\gamma,z\rangle)$ and we recall from \cite{HelR} that
$u$ is \textit{a priori} only $2 \Gamma$-periodic.
Let $\Delta_f$ denote the set of frequencies $\gamma\in \frac{1}{2}\Gamma^*$ such that
$u_\gamma\neq 0$. It is a subset of
\[
\Gamma_{\beta_0}^* = \{\gamma\in \Gamma^*+\frac{\beta_0}{2}:
|\gamma|=\frac{|\beta_0|}{2},\gamma\neq
\pm\frac{\beta_0}{2}\}
\]
which is invariant under $\gamma\mapsto-\gamma$ (see \cite{HelR}). Define also 
$u^+ := \frac{1}{2}(1-J i)u$, the projection on $i$-eigenspace of $J$
(and similarly define $u^-$ for the $(-i)$-eigenspace). Then equation
(\ref{eq:psum}) projects down to the $i$-eigenspace as
\[
\forall \gamma, \quad 0=(i\pi\bar\gamma)^N\sum_{k=0}^N (-\frac{\bar\beta_0}{2\bar\gamma})^kp_k
u^+_\gamma
\]
and similarly for the $(-i)$-eigenspace. One easily shows that $u^+_\gamma$ and $ u^{-}_\gamma$ are both non-zero for $\gamma\in\Delta_f$ so that
\begin{equation}\label{eq:roots}
p(\frac{\bar\beta_0}{2\bar\gamma}) = 0\ \text{whenever}\ \gamma \in \Delta_f.
\end{equation}
Since $|\gamma|=|\beta_0/2|$ these roots of $p(\lambda)$ lie on the unit circle. On top of this,
the $\bar\rho^*$ and $\tau^*$ symmetries of the \pk field equations, together with the uniqueness
result
in lemma \ref{lem:minimal}, tell us that
\begin{eqnarray}\label{eq:xsymmetry}
(\overline{p(\bar\lambda^{-1})} I_4,\lambda\overline{x(\bar\lambda^{-1})}) &  =  &
\overline{p(0)}\lambda^{-N}(p(\lambda) I_4,x(\lambda)), \\
(p(-\lambda) I_4,-iLx(-\lambda)) &  = &  (-1)^N(p(\lambda) I_4,x(\lambda)),\notag
\end{eqnarray}
We note that $\overline{p(0)}=1/p(0)$ and $N$ must be even for $p$ to have minimal degree.
Taken altogether this information allows us to completely determine $p(\lambda)$.
\begin{lem}
The minimal \pk field $(p(\lambda) I_4,x_\lambda)$ has
\begin{equation}\label{eq:p}
p(\lambda) = \prod_{\gamma\in\Delta_f} (\lambda - 2\gamma/\beta_0).
\end{equation}
In particular, $p(\lambda)$ is even, its roots are simple and all lie on the unit circle.
\end{lem}
We will label these roots $s_j$, $j=1,\ldots,N$ and note that
\begin{equation}\label{eq:roots2}
s_j = 2\gamma_j/\beta_0 = \bar\beta_0/2\bar\gamma_j,
\end{equation}
where $\{\gamma_j:j=1,\ldots,N\}=\Delta_f$. It will be convenient for us later to label these so
that,
for $s_{j+N/2}=-s_j$ (i.e., $\gamma_{j+N/2}=-\gamma_j$).
\begin{proof}
Assume that $p(\lambda)$ is given by \eqref{eq:p} and define $x_\lambda$ according to
\eqref{eq:recursion}. This necessarily satisfies \eqref{eq:pkf2}: we must show that it also
satisfies
\eqref{eq:pkf3} (equivalently, that it has the real symmetry in \eqref{eq:xsymmetry}).
But that follows from \eqref{eq:linear} by a straightforward computation.
\end{proof}
\begin{rem}\label{rem:conformal}
A particular consequence of \eqref{eq:recursion} and \eqref{eq:xsymmetry} is that 
$x_N = -\frac{2}{\pi\beta_0}J\bar u = p_0\bar x_1$.
Under the assumption that $f$ is conformal this means $x_1,x_N$ never vanish. Later we will consider
the possibility that $f$ is only weakly conformal: the branch point of $f$ will occur precisely when
the degree of $x_\lambda$ drops.
\end{rem}
Finally, consider the effect of the group of symplectic isometries on the minimal polynomial
Killing field. Since we wish to retain the conditions $f(0)=0$  and $\beta(0)=0$ we are only interested 
in the action of $G_0\subset G$. 
\begin{lem}\label{lem:rotations}
If $f:\C/\Gamma\to\R^4$ has minimal polynomial Killing field $(p(\lambda)I_4,x_\lambda)$ then, for any
$g\in G_0$, $gf$ has minimal polynomial Killing field $(p(\lambda)I_4,gx_\lambda)$.
\end{lem}
The proof follows at once from the construction above.

\section{The spectral data.}\label{sec:spectral}
As with other surface geometries which arise from integrable systems, 
the spectral data for a \HSL torus consists of a complete
algebraic curve, a rational function on that curve and a line bundle over that curve. 
Two features which distinguish this geometry  are that: a) the spectral curve is rational (indeed
reducible), b) the spectral data is not invariant under the action of ambient symmetries (the
symplectic isometries of $\R^4$). Both of these features oblige us to take extra care when
formulating the correspondence between spectral data and \HSL tori.

\subsection{The spectral curve.}
As a general principle (see e.g.\ \cite{McI}) the spectral data should be a geometric realisation
of the algebra $\caK$ of all \pk fields. As a straightforward consequence of lemma \ref{lem:minimal} 
we obtain a complete description of $\caK$, by thinking of it as a module over the algebra $\caB=\C[\lambda^{-1},\lambda]$ in the obvious way.
\begin{lem}\label{lem:K}
Let $(p(\lambda) I_4,x_\lambda)$ be the minimal \pk field of this type for $f$ and let $R\in\End(\C^4)$
be an orthogonal projection commuting with $J$. Then the minimal \pk field for $R$ is
$(\frac{p}{r}R,\frac{1}{r}Rx_\lambda)$ where $r(\lambda)$ is the monic polynomial whose
zeroes are exactly the common zeroes of $p(\lambda)$ and $Rx_\lambda$. Thus $\caK$ is generated over
$\caB$ by the set of all these minimal \pk fields together with the constant \pk field $I$. 
\end{lem}
In particular, since $p(\lambda)$ and $x_\lambda$ do not have
common zeroes, $(p(\lambda) R, Rx_\lambda)$ is the minimal \pk field for $R$ unless $x_\lambda$ 
lies in $\ker(R)$ at some zeroes of $p(\lambda)$. 

Since $\caK$ is non-commutative it is necessary to construct 
the spectral curve by taking a maximal abelian subalgebra $\caA$ of $\caK$, as the
completion by smooth points of the affine curve $\Spec(\caA)$. Given our
symmetries $\caA$ should correspond to a 
choice of real maximal torus $\ft$ in the commutator $\fz^J\subset\End_\R(\R^4)$ of $J$ for which $\Ad
L\cdot\ft = \ft$. Let us first consider the maximal torus $\ft_0\subset \fz^J$ generated by
\[
\begin{pmatrix} I_2 & 0 \\ 0 & 0 \end{pmatrix},\
\begin{pmatrix} J_2 & 0 \\ 0 & 0 \end{pmatrix},\
\begin{pmatrix} 0 & 0 \\ 0 & I_2 \end{pmatrix},\
\begin{pmatrix} 0 & 0 \\ 0 & J_2 \end{pmatrix},
\]
where $I_2$ is the identity, and $J_2$ the standard complex structure, on $\R^2$.
This maximal torus is characterised by its decomposition of $\C^4$ into (Hermitian orthogonal)
invariant lines given by
\begin{equation}\label{eq:V}
V_1 = \C.(\e_1-i\e_2),\ V_2 = \C.(\e_3-i\e_4),\ V_3 = \bar V_1 = LV_2,\ V_4 = \bar V_2 = LV_1.
\end{equation}
Notice that $\ft_0$ is completely determined by $V_1$, which is a complex line in the $i$-eigenspace
$V\subset\C^4$ of $J$.
Every other maximal torus of $\fz^J$ for which $\Ad L\cdot\ft = \ft$ is given by $\Ad g\cdot\ft_0$ for
some $g\in G_0$.

For a fixed $f$ we can define
\[
\caA_f = \{(q(\lambda)R,y_\lambda)+u(\lambda)I\in\caK: R\in\ft_0^\C\}.
\]
By lemma \ref{lem:K} $\caA_f$ is abelian. Further, it is maximal abelian. For if $\xi\in\caK$ 
commutes with every element of
$\caA_f$ then $\xi = (q(\lambda)R,y_\lambda)$ where $R\in\ft_0^\C$ since $\ft_0^\C$ is maximal.

However, quite surprisingly, it is \emph{not} true that $\caA_{gf}\simeq\caA_f$ for every $g\in
G_0$ (recall that $G_0$ is the group of symplectic isometries preserving the base point conditions
$f(0)=0$ and $\beta(0)=0$). This isomorphism is only true generically in $G_0$, as we will show
below.  The generic algebra is isomorphic to an algebra $\caA$ which we will now describe.
For each $k=1,\ldots,4$ let $R_k\in\End_\C(\C^4)$ denote the unitary projection matrix 
corresponding to orthogonal projection onto $V_k$. These span $\ft_0^\C$ and satisfy 
$R_jR_k = \delta_{jk}R_k$ (where $\delta_{jk}$ is the Kronecker delta) and
\begin{equation}\label{eq:R}
R_3 = -LR_2L,\quad R_4 = -LR_1L, \quad R_3 = \bar R_1,\quad R_4 = \bar R_2.
\end{equation}
Now define $\xi_k = (p(\lambda) R_k,R_k x_\lambda)$. Then 
\begin{equation}\label{eq:xi}
\xi_j\xi_k = \delta_{jk}p(\lambda)\xi_k.
\end{equation}
Finally, define
\[
\caA = \C[\xi_1,\xi_2,\xi_3,\xi_4,\lambda I,\lambda^{-1}I]
\]
This is clearly an abelian subalgebra of $\caK$. 
Notice that  any $\eta\in\caA$ can be uniquely written in the form
\[
\eta = \sum_{j=1}^4q_j\xi_j+q_5I,\quad q_j\in\C[\lambda,\lambda^{-1}].
\]
It is easy to see from this that $\eta=0$ if and only if each $q_j=0$ and therefore $\caA$ has no
non-trivial relations other than those in \eqref{eq:xi}. We deduce the following.
\begin{lem}\label{lem:Aisom}
$\caA\simeq \C[Z_1,Z_2,Z_3,Z_4,Z_5,Z_5^{-1}]/\caI$ where $\caI$ is the ideal generated by
$Z_k^2-p(Z_5)Z_k,Z_jZ_k$ for $j,k=1,\ldots,4$ and $j\neq k$.
\end{lem}
The next lemma describes the situation under which $\caA$ is maximal abelian.
\begin{lem}\label{lem:Aequal}
$\caA\subseteq\caA_f$, and $\caA = \caA_f$ if and only if $R_k(x(s_j,0))\neq 0$ for $k=1,\ldots,4$,
$j=1,\ldots,N$.
\end{lem}
\begin{proof}
By definition $\caA\neq \caA_f$ if and only if there exists a minimal polynomial Killing field for 
$R\in\ft_0^\C$,
which we write as $(qR,y)$, not in $\caA$. Since it is minimal there exists $r\in\C[\lambda]$ such that
$(rqR,ry)=(pR,Rx)$. This is equivalent to saying that $p$ and $Rx$ have at least one common zero, i.e., 
there
exists at least one $s_j$ for which $R_k(x(s_j,z))=0$, for some $k=1,\ldots,4$ and for all $z\in\C$.
Further, the common zeroes are independent of $z$, for 
at the zeroes of $p(\lambda)$ \eqref{eq:pkf1} reduces to 
\[
d(\exp(\frac{1}{2}\beta_\lambda J)x_\lambda) = 0,\quad 
\beta_\lambda = \pi(\lambda^{-1}\bar\beta_0z +\lambda\beta_0\bar z),
\]
hence
\begin{equation}\label{eq:xevolution}
x(s_j,z) = \exp(-\frac{1}{2}\beta(s_j,z) J)x(s_j,0),\quad j=1,\ldots,N.
\end{equation}
It follows that
\[
R_k(x(s_j,z)) = e^{\pm i\beta(s_j,z)/2}R_k(x(s_j,0)),
\]
and therefore $\caA\neq\caA_f$ if and only if $R_k(x(s_j,0))= 0$ for some $j,k$.
\end{proof}
Now we can show that generically $\caA=\caA_f$. To be precise, consider the set
\[
\caU_f=\{g\in G_0:R_k(gx(s_j,0))\neq 0\ \forall\ k=1,\ldots,4,\ j=1,\ldots,N\ z\in\C\}.
\]
\begin{lem}\label{lem:U}
The non-empty subset $\caU_f\subset G_0$ is proper, open and $\caA_{gf}\simeq\caA$ if and only if 
$g\in \caU_f$.
\end{lem}
\begin{proof}
Recall that $G_0\subset U(2)\subset SO(4)$ is the commutator of $J$ and $L$. Therefore
$G_0$ preserves both $V=V_1\oplus V_2$ and $\bar V$, the $\pm i$-eigenspaces of $J$.
Further, $G_0$ acts transitively on $\P V$ and the subset of $g\in G_0$ for which
$gV_1\cap V_1=\{0\}$ is an open subset of $G_0$ (hence, the same statements 
are true for the action of $G_0$ on
$\bar V$ and the orbit of $\bar V_1$).  Finally, whenever $|\lambda|=1$ the real
symmetry in \eqref{eq:xsymmetry} implies that the line $\ell = \C.x(\lambda,0)$ has $\bar\ell = \ell$ and 
therefore has non-trivial components in both $V$ and $\bar V$. It follows that for each root
$s_j$ of $p(\lambda)$ there is a proper non-empty open subset $\caU_{kj}\subset G_0$ for which each 
$R_k(gx(s_j,0))$ is non-zero:  $\caU_f$ is the
intersection of the finitely many open subsets $\caU_{kj}$. It follows from lemmas 
\ref{lem:Aisom} and \ref{lem:Aequal} that  $\caA_{gf}\simeq\caA$ if and only if $g\in \caU_f$.
\end{proof}
Given $f$ we define its affine spectral curve to be the affine
scheme $X_0 = \Spec(\caA)$. The natural inclusion of algebras $\caB\to \caA$ is
dual to a finite morphism $\lambda:X_0\to\C^\times$. It follows from \eqref{eq:xi} that $\caA$ is a 
free rank $5$ module over
$\C[\lambda,\lambda^{-1}]$, hence this morphism has degree $5$. We define the spectral curve $X$ to be the
completion of $X_0$ by smooth points.  By lemma \ref{lem:Aisom} $X_0$ is 
biregular to the curve in $\C^5\setminus\{Z_5=0\}$ determined by the equations 
\begin{equation}\label{eq:curve}
Z_j(Z_j-p(Z_5))=0,\ Z_jZ_k=0,\ j,k=1,\ldots,4,\ k\neq j.
\end{equation}
Its normalisation $\varphi:\tilde{X_0}\to X_0$ is dual to the algebra
monomorphism
\[
\varphi^*:\caA\to\tilde\caA;\quad \sum_{j=1}^4q_j\xi_j+q_5I\to (q_1p +q_5,\ldots,q_4p +
q_5,q_5),
\]
where $\tilde\caA$ denotes $\caB^5$ with the direct product structure. 
We obtain the following structure for the spectral curve $X$. 
\begin{prop}
$X$ is a reducible rational curve with five irreducible components
$C_1,\ldots,C_5$ each of which is a smooth rational curve. Any two intersect along the
$0$-dimensional subscheme $\fS\subset X_0$ given by $\fS=\Spec(\caA/\caJ)$ where $\caJ$ is the ideal in
$\caA$ generated by $p,\xi_1,\ldots,\xi_4$. As a divisor, $\fS$ is just the union 
of the $N$ singular points corresponding to the zeroes of $p(\lambda)$. In particular, $X$ has 
arithmetic genus $g=4(N-1)$.
\end{prop}
For simplicity, we will abuse notation by also using $\fS$ to denote the 
set $\lambda(\fS)=\{s_1,\ldots,s_N\}$.

It will be convenient for us to identify the irreducible components as follows. The component 
$C_5$ has affine part given by $Z_j=0$ for $j=1,\ldots,4$, while 
the affine part of $C_j$, for $j\leq 4$, is the curve with equations
\[
Z_j=p(Z_5),\quad Z_k=0,\ k\neq j,\ k\leq 4.
\]
The holomorphic function $\lambda=Z_5$ on $X_0$ extends to a rational function $\lambda:X\to\hat\C$, where
$\hat\C=\C\cup\{\infty\}$. It has
$5$ points over $\lambda=0$, which we will denote by $P_1,\ldots,P_5$, with $P_j$ lying on $C_j$. The 
corresponding points over $\lambda=\infty$ will be denoted $Q_1,\ldots,Q_5$. 

Since $\alpha_\lambda$ possesses the symmetries $\rho^*$ and $\tau^*$ these act as involutions on
$\caA$. Therefore they induce involutions on $X$, which we will call $\rho$ and
$\tau$, for which $\overline{h\circ\rho} = \bar\rho^* h$ and $h\circ\tau = \tau^*h$ for any
$h\in\C[X_0]$. It is straightforward to establish the following characterisation of these involutions.
\begin{prop}\label{prop:rhotau}
The involutions $\rho$ and $\tau$ on $X$ are, respectively, anti-holomorphic and holomorphic, and act as follows.
\begin{eqnarray*}
(\overline{Z_1\circ\rho},\overline{Z_2\circ\rho},\overline{Z_3\circ\rho},\overline{Z_4\circ\rho},
\overline{Z_5\circ\rho})&  = &
(\frac{Z_3}{Z_5^Np(0)},\frac{Z_4}{Z_5^Np(0)},\frac{Z_1}{Z_5^Np(0)},\frac{Z_2}{Z_5^Np(0)},
\frac{1}{Z_5}),\\
\tau(Z_1,Z_2,Z_3,Z_4,Z_5) &= &(Z_4,Z_3,Z_2,Z_1,-Z_5).
\end{eqnarray*}
\end{prop}
Finally, let us explain the difference between $\Spec(\caA_f)$ and $\Spec(\caA)$ when
$\caA\neq\caA_f$. For each $k=1,\ldots,4$ define
\begin{equation}\label{eq:fS_k}
\fS_k = \{s_j:R_k(x(s_j,0))\neq 0\}.
\end{equation}
Note that every $s_j$ belongs to at least one $\fS_k$ since $p(\lambda)$ and $x(\lambda,0)$ have no 
common zeroes. Now define
\[
p_k(\lambda) = \prod_{s_j\in \fS_k}(\lambda - s_j).
\]
As a consequence of the symmetries \eqref{eq:xsymmetry} and \eqref{eq:R} we see that 
$\fS_1=\fS_3$, $\fS_2=\fS_4$
and $\fS_2=-\fS_1$. Therefore $p_3=p_1$, $p_4=p_2$ and $p_2(\lambda) = \pm p_1(-\lambda)$. 
It follows that $\caA_f = \C[\eta_1,\ldots,\eta_4,\lambda I,\lambda^{-1}I]$ where
\[
\eta_k = \frac{p_k}{p}\xi_k = (p_kR_k,\frac{p_k}{p}R_kx_\lambda)
\]
and therefore $\Spec(\caA_f)$ is biregular to the affine curve in $\C^5\setminus\{Z_5=0\}$ with
equations
\[
Z_j(Z_j-p_j(Z_5))=0,\ Z_jZ_k=0,\ j,k=1,\ldots,4,\ k\neq j.
\]
We conclude from this the following structure.
\begin{lem}\label{lem:X_f}
Let $X_f$ denote the completion by smooth points of $\Spec(\caA_f)$. 
It has five irreducible components $C_j^f$, $j=1,\ldots,5$, each with $C_j\simeq\hat\C$.
All five components intersect along $\fS_1\cap \fS_2$, with further intersection 
relations
\[
C_1^f\cap C_3^f=\fS_1,\ C_2^f\cap C_4^f=\fS_2,\ C_k^f\cap C_5^f=\fS_k,\ k=1,\ldots,4.
\]
In particular, $X_f$ has arithmetic genus $4(N_1-1)$, where $N_1 = \# \fS_1$. 
\end{lem}
The picture one should have in mind is that
the natural inclusion $\caA\hookrightarrow\caA_f$ is dual to a
finite morphism $\pi_f:X_f\to X$ which realises $X_f$ as the desingularisation of 
$X$ obtained by pulling apart the components $C_1,\ldots,C_5$ so that only the above 
intersection relations remain and the symmetries $\rho$ and $\tau$ persist. 
Notice that, because of the symmetries $\rho$ and $\tau$, the structure of
$X_f$ is completely determined by knowing how $C_1^f$ and $C_5^f$ intersect.
\begin{rem}
In case this seems an overly elaborate way of obtaining the spectral curve, consider the alternatives.
Taking the characteristic polynomial of $\xi_\lambda = (p(\lambda)I_4,x_\lambda)$ we obtain
an unreduced planar curve with equation $\mu(\mu-p(\lambda))^4=0$. A slightly more sophisticated
approach is to consider the curve of eigenlines of $\xi_\zeta$, but its eigenlines are generated by
$\e_1,\ldots,\e_4$, which are constant, and the vector 
\[
(x_1(\lambda,z),\ldots,x_4(\lambda,z),-p(\lambda)),\quad x = \sum_{j=1}^4x_j\e_j.
\]
This leads to a disconnected union of Riemann spheres. Neither of these approaches permits
the vector $x(\lambda,z)$ to be encoded in the spectral data using a line bundle, or sheaf, over the
spectral curve. The next section show how our spectral data achieves this.
\end{rem}
\subsection{Line bundles over the spectral curve.}
In actuality $\caA$ is a family of  algebras. By evaluating every \pk field at a point $z$ we
construct an algebra $\caA(z)$: these are all isomorphic under the map
\begin{equation}\label{eq:conj}
\caA(0)\to\caA(z);\ \xi_\lambda(0)\mapsto\Ad U_\lambda(z)^{-1} \xi_\lambda(0),
\end{equation}
where $U_\lambda$ is the extended frame, i.e., $U_\lambda^{-1}dU_\lambda = \alpha_\lambda$ and
$U_\lambda(0)=I$.
The vector space $\caM=\caB\otimes\C^5$ is an $\caA(z)$-module for each $z$. We will will show that it
determines a line bundle over each irreducible component and then explain how these fit together
over $X$. Before we do this, we need to review the moduli space $\Pic(X)$ of line bundles over $X$. 

Let $\tilde X$ be the normalisation of $X$: it is the disjoint union of 
the smooth rational curves $C_1,\ldots,C_5$. A line bundle over $X$ can be thought of as a
divisor equivalence class $[D]$ for a divisor $D$ of smooth points on $X$. Recall (from e.g.,
\cite{Ser}) that the equivalence is characterised by the property of being trivial if and only if $D$ is
the divisor of a rational function on $X$. Thus $\Pic(X)$ is
isomorphic to the group of these divisor equivalence classes. Since $X$ is reducible such a class
cannot be assigned a single integer for its degree, but rather
\[
\deg:\Pic(X)\to\Z^5;\quad \deg([D]) = (\deg([D\cap C_1],\ldots,\deg([D\cap C_5])).
\]
In particular, the Jacobi variety $\Jac(X)$ is the subgroup of line bundles whose degree on each 
component is zero. The structure of $\Jac(X)$ is given by the next lemma, whose proof is
straightforward. 
\begin{lem}\label{lem:Jac(X)} 
Let $\Div_0(X)$ denote the set of all divisors with $\deg(D)=(0,\ldots,0)$.  
For each $D\in\Div_0(X)$ there is a rational function $f$ on $\tilde X$, unique up to
scale, with divisor $D$.
The map $\Jac(X)\to (\Ct)^{4(N-1)}$ which assigns to $[D]$ the coordinates 
\begin{equation}\label{eq:coords}
t_{kj} = \frac{f_k(s_j)}{f_k(s_N)}\frac{f_5(s_N)}{f_5(s_j)},\quad 
1\leq k\leq 4,\ 1\leq j\leq N-1,
\end{equation}
in which $f_k=f|_{C_k}$, is an isomorphism of linear algebraic groups.
\end{lem}
\begin{lem}\label{lem:caL}
The $\caA(z)$-module $\caM$ determines a rank $1$ sheaf $\caL_f(z)$ over $X$. Either:
\begin{enumerate}
\item $\caA(z) = \caA_f(z)$ and $\caL_f(z)$ is a line bundle of degree $(N,N,N,N,0)$, or,
\item $\caA(z)\neq\caA_f(z)$ and $\caL_f(z)$ is the direct image of a line bundle over
$X_f$ of degree $(N_1,N_1,N_1,N_1,0)$.
\end{enumerate}
In either case $\caL_f$ possesses the symmetries
\begin{equation}\label{eq:caLsymmetry}
\overline{\rho^*\caL_f} \simeq \caL_f(P_5-Q_5),\quad \tau^*\caL_f\simeq \caL_f.
\end{equation}
\end{lem}
For the purpose of the proof, and for subsequent use, it will be convenient to work 
with the Hermitian orthonormal basis $v_1,\ldots,v_4$ for $\C^4$, thought of as the 
first four dimensions in $\C^5$, with $v_j\in V_j$ defined by
\begin{equation}\label{eq:basis}
v_1 = \frac{1}{\sqrt{2}}(\e_1-i\e_2),\ v_2 = \frac{1}{\sqrt{2}}(\e_3-i\e_4),\ 
v_3 = \bar v_1,\ v_4 = \bar v_2. 
\end{equation}
We note that $v_2=L\bar v_1,v_4=Lv_1$.
In such a basis we write $x_\lambda(z) = \sum_{j=1}^4\chi_j(\lambda,z)v_j$, so that
$R_kx_\lambda = \chi_k(\lambda)v_k$. We also define  $v_5 =
(0,0,0,0,1)$.
\begin{proof}
First assume $\caA=\caA_f$ (for simplicity we drop the explicit dependence of the notation
on $z$). Let $\caI_j\subset\caA$ denote the prime ideal corresponding 
to the component $C_j\cap X_0$. In terms of generators we have
\begin{eqnarray*}
\caI_j& = &\langle\xi_j-p I,\xi_k:k\neq j\rangle,\quad j\leq 4\\
\caI_5 & = & \langle \xi_1,\ldots,\xi_4\rangle
\end{eqnarray*}
Now define $\caA_j=\caA/\caI_j$ (this is the coordinate ring for $C_j\cap X_0$) and
$\caM_j = \caM/\caI_j\caM$. It is easy to see that $\caA_j\simeq\caB$ for each 
$j=1,\ldots,5$ and the points over $\lambda=0,\infty$ on $C_j$ correspond to the two 
gradings carried by $\caA_j$, namely, the degree in $\lambda^{-1}$ and $\lambda$, respectively.
For $j=1,\ldots,4$ we have
\[
\caI_j\caM = \caB\langle p v_k,\chi_kv_k,\chi_jv_j-p v_5:k\neq j\rangle.
\]
But $p(\lambda)$ and $\chi_k(\lambda)$ have no common zeroes so
the ideal in $\caB$ generated by $p$ and $\chi_k$ is $\caB$ itself. Therefore
\[
\caI_j\caM = \caB\langle v_k,\chi_jv_j-p v_5:k\neq j\rangle.
\] 
It follows that
\[
\caM_j\simeq \caB\langle v_j,v_5\rangle/\langle \chi_jv_j-p v_5\rangle.
\]
This is clearly a rank one module over $\caA_j$ and torsion free since $p$ and $\chi_j$ have
no common zeroes. This determines a line bundle $\caL_j$ over $C_j$ with sections
$\sigma_j,\sigma_5$, corresponding to $v_j,v_5$, satisfying $\chi_j\sigma_j-p\sigma_5=0$. 
Now
\[
\caI_5\caM = \caB\langle v_1,\ldots,v_4\rangle
\]
and therefore $\caM_5\simeq \caB\langle v_5\rangle$, which is clearly a rank one torsion 
free module over $\caA_5$. It determines a trivial line bundle $\caL_5$ with nowhere
vanishing global section $\sigma_5$. 

Therefore as an $\caA$-module $\caM$ determines a sheaf $\caL_f$ over $X$ for which 
$\caL_f|C_j=\caL_j$. It has globally holomorphic section $\sigma_5$ which 
vanishes exactly at the zeroes of $\chi_j(\lambda,z)$: 
this gives a degree $N$ divisor $D_j(z)$ on $C_j$, so $\deg(\caL_j)=N$. Since $\caL_5$ is
trivial it has degree $0$. Hence $\sigma_5$ has divisor
\begin{equation}\label{eq:D}
D(z) = D_1(z)+\ldots+ D_4(z).
\end{equation} 
Since this divisor includes no singular points 
$\caL(z)$ is invertible (i.e., a line bundle) with $\caL(z)\simeq\caO_X(D(z))$.

Now consider the case $\caA\neq\caA_f$.  A simple adaptation of the arguments above shows
that the $\caA_f$-module $\caM$ determines a line bundle over $X_f$ of degree
$(N_1,N_1,N_1,N_1,0)$. Therefore its direct image $\caL_f$ corresponds to $\caM$ as an
$\caA$-module.

Finally, from the symmetries \eqref{eq:xsymmetry} and \eqref{eq:basis} we obtain
\begin{eqnarray}\label{eq:chi} 
\chi_3(\lambda) = \lambda^{N+1}p(0)\overline{\chi_1(\bar\lambda^{-1})}, &
\chi_4(\lambda) = \lambda^{N+1}p(0)\overline{\chi_2(\bar\lambda^{-1}}), \\
\chi_4(\lambda) = -i\chi_1(-\lambda), & \chi_3(\lambda) = i\chi_2(-\lambda).\notag
\end{eqnarray}
The first expression gives the equation of divisors
\[
D_3-NQ_3 = (N+1)(P_3-Q_3)+\rho^*(D_1-NQ_1).
\]
It follows that $D_3+Q_3-P_3=\rho^*D_1$. We obtain a similar equation relating $D_4$ and $\rho^*D_2$,
and deduce
\[
\rho^* D = D + \sum_{j=1}^4(Q_j-P_j)\sim D+P_5-Q_5,
\]
using the fact that $\sum_{j=1}^5(P_j-Q_j)$ is the divisor of $\lambda$. Similarly $\tau^*D=D$. 
\end{proof}
\begin{rem}
In this lemma we are implicitly assuming that the 
degree of $\chi_j(\lambda,z)$ is exactly $N$ for each $z\in\C$. This means that 
$D_j(z)$ actually lies on $C_j\setminus\{Q_j\}$. But it is possible for the degree of $\chi_j$ to drop
(e.g., at branch points of $f$: see below). We shall 
keep $\deg(D_j(z))=N$ by allowing it to include the points at $\infty$.  
\end{rem}
It is an inevitable consequence of lemma \ref{lem:U}
that the group $G_0$ acts non-trivially on the spectral data.
Indeed, the $G_0$-orbit of $f$ takes $\caL_f$ outside
the Jacobi variety $\Jac(X)$ into a compactification which includes non-invertible sheaves,  
since there is at least one point in this orbit 
at which $\chi_k$ has a zero at some $s_j$.  We will
explain precisely what happens when we describe the moduli space of HSL tori in the next section. 
Meanwhile let us observe that the previous proof provides a convenient characterisation of the sheaf 
$\caL_f$, since it shows that $\caL_f$ is entirely determined by the divisor of zeroes of $\chi_1$, 
even when these are allowed to include singularities or points at infinity.
\begin{cor}\label{cor:divisor}
Up to isomorphism $\caL_f$ is determined by the positive divisor $E=D_1-P_1$ of degree
$N-1$ on $C_1$, where $D_1$ is the divisor of zeroes of $\chi_1$. The map $f$ has a branch point at
$z=0$ precisely when $E$ has the form $E=P_1+Q_1+E'$ for some positive divisor $E'$ of degree $N-3$.
\end{cor}
\begin{proof}
From \eqref{eq:pkf2} and \eqref{eq:pkf3} we know $x_\lambda$, and hence $\chi_1$, has a zero at $\lambda=0$.
From remark \ref{rem:conformal} we know $f$ has a branch point at $z=0$ precisely when the degree of
$x_\lambda$ drops. Using the symmetries \eqref{eq:chi} we see that this occurs precisely when $\chi_1$
has an additional zero at $P_1$ and degree less then $N-1$, i.e., a zero at $Q_1$.
\end{proof}
We finish this section by describing how the sheaf $\caL_f(z)$ moves with base point translation. 
Since we started with a Lax
equation we expect this motion to be linear, and that is indeed the case.
\begin{prop}
The map $\ell:\C/\Gamma\to\Jac(X)$ given by $\ell=\caL(z)\otimes\caL(0)^{-1}$ has image in
the real analytic subgroup $J_R\subset\Jac(X)$ consisting of all line bundles possessing the
symmetries
\begin{equation}\label{eq:Lsymmetry}
\overline{\rho^*\ell} \simeq \ell,\quad \tau^*\ell\simeq \ell.
\end{equation}
Further, $\ell$ is a homomorphism of real groups characterised by the equation
\begin{equation}\label{eq:Ltangent}
\frac{\partial \ell}{\partial z}|_{z=0} = 
\frac{i\pi\bar\beta_0}{2}\frac{\partial\caA}{\partial\lambda}|_{\lambda=0}
\end{equation}
where $\caA:\C\to\Jac(X)$ is defined by
\begin{equation}\label{eq:A}
\caA(\lambda) = \caO(\sum_{j=1}^4(P_j(\lambda)-P_j)),
\end{equation}
with $P_j(\lambda)$ the point on $C_j$ given by $Z_j=\lambda$ for $j=1,2$ but $Z_j=-\lambda$ for
$j=3,4$. 
\end{prop}
\begin{proof}
The isomorphisms in \eqref{eq:Lsymmetry} are a direct consequence of \eqref{eq:caLsymmetry}.
To prove \eqref{eq:Ltangent} we use the coordinates 
\eqref{eq:coords} on $\Jac(X)$, in which $L(z)$ has coordinates,
\[
t_{kj}(z) = \frac{\chi_k(s_j,z)}{\chi_k(s_j,0)}\frac{\chi_k(s_N,0)}{\chi_k(s_N,z)},\quad
1\leq k\leq 4,\ 1\leq j\leq N-1.
\]
Here we are using the fact that $D(z)-D(0)$ is the divisor for the rational function $\chi$ on $\tilde
X$ whose restriction to $\tilde C_k$ is $\chi_k$ for $1\leq k\leq 4$. 
Now using \eqref{eq:xevolution} we compute
\begin{equation}\label{eq:dt/dz}
\frac{\partial t_{kj}}{\partial z}|_{z=0} = \epsilon_k\frac{-i\pi\bar\beta_0}{2}(s_j^{-1}-s_N^{-1}),
\end{equation}
where $\epsilon_k$ equals $+1$ for $k=1,2$ and $-1$ for $k=3,4$. On the other hand
we compute, in analogy with the computation in the appendix, the coordinates of $\caA(\lambda)$
to be
\[
a_{kj}(\lambda) = \frac{\epsilon_k\lambda-s_j}{\epsilon_k\lambda-s_N}\frac{s_N}{s_j}, \quad
1\leq k\leq 4,\ 1\leq j\leq N-1.
\]
Therefore
\[
\frac{\partial a_{kj}}{\partial\lambda}|_{\lambda=0} = \epsilon_k(s_N^{-1}-s_j^{-1}).
\]
Equation \eqref{eq:Ltangent} follows.
\end{proof}
The real subroup $J_R$ acts on the set of spectral data by tensor product $\caL\mapsto\caL\otimes L$,
for $L\in J_R$. This action gives  what are usually called the ``higher flows'' in integrable systems
language. 
\begin{lem}\label{lem:J_R}
The real subgroup $J_R\subset\Jac(X)$ is isomorphic, as a real group, to 
$\Ct)^{N-1}$. In the coordinates \eqref{eq:coords} on $\Jac(X)$ it corresponds to the subgroup
\[
\{t_{kj}\in\Ct: t_{1j}=\frac{t_j}{t_N},\ t_{2j} =\frac{\bar t_{j+n}}{\bar t_{N/2}},\ 
t_{3j}=\frac{\bar t_j}{\bar t_N},\ t_{4j}=\frac{t_{j+N}}{t_{N/2}},\ \exists t_j\in\Ct,\ j\in\Z_N\}.
\]
\end{lem}
This is straightforward to prove: the action of $\Jac(X)$ equates to the action of $(\Ct)^{4(N-1)}$ as
the group of diagonal matrices in $(SL_N)^4$ on the $4N$ vector $\chi_k(s_j)$. The subgroup $J_R$ is
that which preserves the symmetries \eqref{eq:chi}.

\section{The moduli spaces of spectral data and of \HSL tori.}
\subsection{Reconstruction from the spectral data.}
We begin by showing that the spectral data reconstructs the HSL torus $f$, up to dilations
and base point preserving
isometries which preserve the maximal torus $\ft_0$. These isometries comprise
\[
T_0 = \{g\in G_0:\Ad g\cdot\ft_0=\ft_0\}\simeq S^1.
\]
Fix a Maslov class $\beta_0\in\Gamma^*$ for which $\Gamma^*_{\beta_0}$ is non-empty, and choose a
non-empty subset $\Delta\subseteq\Gamma^*_{\beta_0}$ which is closed under $\gamma\mapsto-\gamma$. 
Now define
\begin{equation}\label{eq:pDelta}
p_\Delta(\lambda) = \prod_{\gamma\in\Delta}(\lambda - 2\gamma/\beta_0)
\end{equation}
and let $\fS_\Delta = \{2\gamma/\beta_0:\gamma\in\Delta\}$ (whose elements we label $s_1,\ldots,s_N$ as
above). Let $X_\Delta$ denote the completion by smooth points of the curve given by using $p_\Delta$ in
equations
\eqref{eq:curve}, and let $\lambda$ be the rational function on it given by $Z_5$. Let $\caE_\Delta$ denote
the set of all positive divisors $E$ of degree $N-1$ on the smooth component $C_1$ of $X_\Delta$ and
with the property that if $s_j$ lies on $E$ then $-s_j$ does not. Let $\overline{\Pic(X)}$ be the set
of rank 1 coherent sheaves over $X$. We can define a map
\begin{equation}\label{eq:L(E)}
\caL:\caE_\Delta\to \overline{\Pic(X_\Delta)},\quad
E\mapsto \caO(E' + \rho^*E' + \tau^*E'+(\rho\tau)^*E' + \sum_{k=1}^4P_k).
\end{equation}
where 
\[
E'=E-\sum_{s_j\in\mathrm{Supp}(E)}s_j
\]
and whenever $E'\neq E$ we mean $\caL(E)$ is the direct image of that line bundle over the partial
desingularisation of $X_\Delta$ obtained in the manner of lemma \ref{lem:X_f}.
\begin{lem}\label{lem:E}
The map \eqref{eq:L(E)} is injective.
\end{lem}
\begin{proof}
Suppose $\caL(E)=\caL(F)$ for $E,F\in\caE_\Delta$. Then there is a rational function $h$ on $X$ for
which $h.\caL(E)_P=\caL(F)_P$ at every stalk. So both must be invertible on $X$ if either is, in which
case $h$
has divisor $E-F$. Therefore $h$ is constant on the component $C_5$ and hence constant globally since
it has degree $N$ on the other components. When $\caL(E)$ is non-invertible, it and $\caL(F)$ must both
come by direct image of a line bundle over the same curve, and these are therefore equivalent line 
bundles upstairs. Again $h$ is constant on $C_5$ and on the other components has degree equal 
to the number of intersection points, hence $h$ is constant globally.
\end{proof}
\begin{prop}\label{prop:reconstruct}
To each triple $(X_\Delta,\lambda,\caL(E))$, with $E\in\caE_\Delta$, there corresponds a 
based map $f:\C/\Gamma\to\R^4$, determined up to dilations and the action of $T_0$, which is a 
weakly conformal HSL immersion with Maslov class $\beta_0$ for which the triple is its spectral data.
\end{prop}
\begin{proof}
Given $E$ there is a unique, up to complex scaling, function $\chi_1(\lambda,0)$ whose divisor is
$E+P_1-N.\infty$. Using the relations \eqref{eq:chi} this gives 
$x(\lambda,0)=\sum\chi_jv_j$. Now define $x(s_j,z)$ using \eqref{eq:xevolution}.
Then $x(\lambda,z)$ is uniquely determined by all $x(s_j,z)$ by linear algebra since
$p(\lambda)$ has $N$ zeroes and $x(\lambda,z)$ has at most $N$ 
non-trivial coefficients $x_1,\ldots,x_N$ as a polynomial in $\lambda$. Finally, set 
\begin{equation}\label{eq:fevolution}
f(z) = \frac{1}{p(1)}(\exp(\frac{1}{2}\beta J)x(1,z) - x(1,0)).
\end{equation}
This is a based map into $\R^4$ and
it is easy to check that effect of the scaling is to dilate and rotate $f$, with rotations from $T_0$.

It remains to show that this map is a HSL immersion.
Since by definition we have $\beta_{z\bar z}=0$ it suffices to show that $Sf_z=if_z$, where
$S=e^{J\beta/2}Le^{-J\beta/2}$ (cf.\ remark \ref{rem:S-hol}). 
Differentiating \eqref{eq:fevolution} gives
\[
d(\exp(\frac{1}{2}\beta J) x(1,z)) = p(1)df,
\]
which is \eqref{eq:pkf1} at $\lambda = 1$. Since $x(\lambda,z)$ also satisfies
\eqref{eq:xevolution} equation \eqref{eq:pkf1} must hold (with $X_\lambda=p_\Delta(\lambda)I$)
at the $N+1$ values $\lambda = 1,s_1,\ldots,s_N$,
and therefore holds for all $\lambda$ since both sides are polynomials in $\lambda$ with at most $N+1$ 
non-trivial terms. Consequently \eqref{eq:pkf2} holds. Since $p_\Delta(\lambda)$ is an even polynomial,
its coefficient $p_j\neq 0$ only if $j$ is even. When $p_j$ is non-zero we have
\begin{equation}\label{eq:pkf4}
e^{-J\beta/2}f_z = \frac{1}{p_j}(\partial x_j/\partial z + \frac{\pi}{2}\bar\beta_0Jx_{j+1}).
\end{equation}
Now we use the symmetry $Lx(-\lambda,z) = ix(\lambda,z)$ to deduce that for $j$ even $Lx_j=ix_j$ and
$Lx_{j+1}=-ix_{j+1}$. So when $L$ is applied to both sides of \eqref{eq:pkf4} we obtain
\[
Le^{-J\beta/2}f_z = ie^{-J\beta/2}f_z,
\]
using $LJ=-JL$. The real symmetry in \eqref{eq:chi} ensures that $f$ takes values in $\R^4$ and the 
$\Gamma$-periodicity is ensured by the fact that the Fourier frequencies of $f$ come from $\Gamma^*$.

Finally, the construction above is clearly the reconstruction from spectral data since it works via the
\pk field $(p_\Delta(\lambda)I,x_\lambda)$.
\end{proof}

\subsection{The moduli space of spectral data.}\label{subsec:spectral} 
From now on we fix a torus $\C/\Gamma$ and a Maslov class $\beta_0\in\Gamma^*$. 
We fix $N=\#\Gamma^*_{\beta_0}$ and set
\[
p(\lambda) = \prod_{j=1}^N(\lambda-s_j),
\]
so that $X$ is now the most singular curve possible to obtain a HSL torus with Maslov class $\beta_0$.
From this point of view, $(X,\lambda)$ is completely determined by $(\Gamma,\beta_0)$.
As before, we identify the singular set $\fS$ with $\{s_j:j=1,\ldots,N\}$.
For every non-empty $\Delta\subset\Gamma^*_{\beta_0}$ of the type used in the previous section 
$X_\Delta$ is a partial desingularisation of $X$, obtained by pulling apart the intersections at
$\fS\setminus\fS_\Delta$. By proposition \ref{prop:reconstruct} the moduli space of spectral data for
HSL tori with Maslov class $\beta_0$ is $\caS(\Gamma,\beta_0)=\cup_\Delta\caE_\Delta$. Also by that
proposition, we can think of $\caS(\Gamma,\beta_0)$ as the moduli space of based, weakly conformal HSL
immersions $f:\C/\Gamma\to\R^4$, with Maslov class $\beta_0$, up to dilations and isometries from $T_0$. 
We will show that this is naturally isomorphic to $\CP^{N-1}$ and provide cooordinates in which the action
of $J_\R\subset\Jac(X)$ corresponds to the natural action on Fourier components.

First, let $\Div_{N-1}(C_1)$ denote the set of all positive divisors of degree $N-1$ on $C_1$. For each
$\Delta$ the map
\begin{equation}
\caE_\Delta\to\Div_{N-1}(C_1);\quad E\mapsto E+\sum_{s_j\not\in\fS_\Delta}s_j,
\end{equation}
is clearly injective. In fact it gives a bijection from $\caS(\Gamma,\beta_0)$ to $\Div_{N-1}(C_1)$,
since $E\in\caE_\Delta$ is not permitted to contain both $s_j$ and $-s_j$ for any $s_j$. 
On the other hand, $\Div_{N-1}(C_1)$ is clearly isomorphic to the projective space 
\[
\P F_{N-1}=\P \{h(\lambda)\in\C[\lambda]:\deg(h)\leq N-1\}
\]
since every polynomial is determined up to scale by its divisor of zeroes. Now we require a result
(which must surely be classical) whose proof is given in appendix \ref{app:theta}.
\begin{prop}\label{prop:varphi}
The map 
\begin{equation}\label{eq:iota}
\iota:\P F_{N-1}\to \CP^{N-1};\quad [h]\mapsto [h(s_1),\ldots,h(s_N)],
\end{equation}
is an isomorphism whose inverse is given by $[\bt]\mapsto [\varphi(\lambda,\bt)]$ where
\begin{equation}\label{eq:varphi}
\varphi(\lambda,\bt) = 
(\prod_{j=1}^N(\lambda-s_j))
\det\begin{pmatrix}
1 & s_1 & \ldots & s_1^{N-2} & t_1/(\lambda-s_1)\\
\vdots & \vdots & &\vdots& \vdots\\ 
1 & s_N&\ldots&s_N^{N-2} & t_N/(\lambda-s_N)
\end{pmatrix},
\end{equation}
whenever $[t_1,\ldots,t_N]=[\bt]$.
The polynomial $\varphi(\lambda,\bt)$ has degree less than $N-1$ if and only if $[\bt]$ 
lies on the hyperplane
\begin{equation}\label{eq:thetainfty}
\Theta_\infty = \{[\bt]\in\CP^{N-1}: \sum_{j=1}^N (\prod_{m\neq j}\frac{1}{s_m-s_j})t_j=0\},
\end{equation}
and has a zero at $\lambda=0$ precisely when $[\bt]$ lies on the hyperplane
\begin{equation}\label{eq:thetazero}
\Theta_0 = \{[\bt]\in\CP^{N-1}: \sum_{j=1}^N (\prod_{m\neq j}\frac{1}{s_m-s_j})\frac{t_j}{s_j}=0\}.
\end{equation}
\end{prop}
Consequently we have fixed an isomorphism $\caS(\Gamma,\beta_0)\simeq\CP^{N-1}$, which we think of as
assigning (homogeneous) coordinates to each based HSL immersion $f:\C/\Gamma\to\R^4$ with Maslov class
$\beta_0$, and the proposition shows how to invert this. Before we describe the inversion, note
that as a corollary to lemma \ref{lem:J_R} the action of the real subgroup $J_R\subset\Jac(X)$ on 
$\caS(\Gamma,\beta_0)$ corresponds in these coordinates to the standard action of $(\Ct)^N/\Ct$ on
$\CP^{N-1}$ by diagonal action on homogeneous coordinates precomposed with the isomorphism
\begin{equation}\label{eq:higherflows}
(\Ct)^{N-1}\to(\Ct)^N/\Ct;\quad (a_1,\ldots,a_{N-1})\mapsto [a_1,\ldots,a_{N-1},1]. 
\end{equation}
In particular, it follows from \eqref{eq:dt/dz} that the action of base point translation on
$S(\Gamma,\beta_0)$ corresponds 
to the action of $\C/\Gamma$ via the real homomorphism
\begin{equation}\label{eq:ell}
\ell:\C/\Gamma\to (\Ct)^N/\Ct;\quad z+\Gamma\mapsto [e_{-\gamma_1}(z),\ldots,e_{-\gamma_N}(z)],
\end{equation}
where we recall that $e_\gamma(z) = \exp(2\pi i\langle\gamma,z\rangle)$. 

Now we describe the explicit inversion of the isomorphism from $\caS(\Gamma,\beta_0)$
to $\CP^{N-1}$ which follows from proposition \ref{prop:varphi}. First, by combining 
\eqref{eq:fevolution} with \eqref{eq:basis} and \eqref{eq:chi} we have the formula
\begin{multline}\label{eq:fconstruction}
f(z) = 2\Re[\frac{1}{p(1)}(\chi_1(1,z)e_{\beta_0/2} - \chi_1(1,0))v_1]\\
+2\Im[\frac{1}{p(1)}(\chi_1(-1,z)e_{-\beta_0/2}-\chi_1(-1,0))Lv_1],
\end{multline}
Thus, as a consequence of corollary \ref{cor:divisor} and 
proposition \ref{prop:varphi}, we achieve the coordinate inversion and characterise the locus of branch
points of $f$.
\begin{prop}\label{prop:bt}
The isomorphism $\caS(\Gamma,\beta_0)\to\CP^{N-1}$ given above is inverted by assigning to 
$[\bt]=[t_j]\in\CP^{N-1}$ the map $f$ given by \eqref{eq:fconstruction}, in which 
\[
\chi_1(\lambda,z)=\lambda\varphi(\lambda,\bt(z)),\quad t_j(z) = t_je_{-\gamma_j}(z).
\]
Further, $f$ has a branch point at $f(z)$ precisely when $[\bt(z)]$ lies on the intersection of
hyperplanes $\Theta_0\cap\Theta_\infty$.
\end{prop}
Later we will require the explicit expression for $f$, which follows from
\begin{eqnarray}
\chi_1(1,z) & = & \sum_{j=1}^N (-1)^{j-1}\Delta_jt_je_{-\gamma_j},\label{eq:chi1} \\
\chi_1(-1,z) & = & \sum_{j=1}^N(-1)^{j-1}\tilde\Delta_jt_je_{-\gamma_j}, \label{eq:chi-1}
\end{eqnarray}
where 
\[
\Delta_j  =  \prod_{\substack{0\leq k<\ell\leq N \\ k,\ell\neq j}}(s_\ell-s_k), \quad
\tilde\Delta_j =\prod_{\substack{1\leq k<\ell\leq N+1 \\ k,\ell\neq j}}(s_\ell-s_k) 
\]
with $s_0=1$ and $s_{N+1}=-1$. 

\begin{rem}\label{rem:branchpts}
Proposition \ref{prop:bt} above allows to make some general conclusions about the existence of branch
points, by exploiting the projective geometry. For simplicity consider the case where
$\beta_0/2\in\Gamma^*$, so that $\Gamma^*_{\beta_0}\subset\Gamma^*$. This is the case referred to in
\cite{HelR} as ``truly periodic'', in that the immersion does not factor through a covered torus.
In that case the real homomorphism $\ell$ in \eqref{eq:ell} has a
complexification obtained in the following way. Fix
generators $\tau_1,\tau_2$ for $\Gamma$ and identify $\C/\Gamma$ with $S^1\times S^1$ via
\[
\C/\Gamma\to S^1\times S^1\subset\Ct\times\Ct;\quad 
u\tau_1+v\tau_2+\Gamma\mapsto (e^{2\pi i u},e^{2\pi i v}),\ u,v\in\R.
\]
The complexification of $\ell$ is 
\begin{equation}\label{eq:ellc}
\ell^\C:\Ct\times\Ct\to(\Ct)^N/\Ct;\quad (a,b)\mapsto [a^{l_1}b^{m_1},\ldots,a^{l_N}b^{m_N}],
\end{equation}
where $l_j= -\langle\gamma_j,\tau_1\rangle$ and $m_j=-\langle\gamma_j,\tau_2\rangle$ are integers. 
This gives a
homomorphism of algebraic groups. Through every point $[\bt]\in\CP^{N-1}$ this puts a
$\Ct\times\Ct$-orbit $\caO_\bt$.  This extends to provide an algebraic map 
\begin{eqnarray*}
\ell_\bt^\C:\CP^1\times\CP^1&\to&\CP^{N-1},\\ 
([a_0,a_1],[b_0,b_1])&\mapsto &
[(\frac{a_1}{a_0})^{l_1}(\frac{b_1}{b_0})^{m_1}t_1,\ldots,(\frac{a_1}{a_0})^{l_N}(\frac{b_1}{b_0})^{m_N}t_N]
\end{eqnarray*}
whose image is the Zariski closure $\bar\caO_\bt$ of the orbit. We can think of $\bar\caO_\bt$ as foliated 
by $\C/\Gamma$-orbits. This gives a $(\R^+)^2$ family, although some will be degenerate orbits. 

For generic $[\bt]$ the pull-back along $\ell^\C_\bt$ of the hyperplanes $\Theta_\infty,\Theta_0$ 
determines curves $\Sigma_\infty,\Sigma_0$ in 
$\CP^1\times\CP^1$. Since $\Theta_\infty$ and $\Theta_0$ are clearly transverse so are these curves,
therefore for generic $[\bt]$ the $\C/\Gamma$-orbit of $[\bt]$ will miss these points, and the
corresponding HSL map will be an unbranched immersion.

For any branched immersion lying in $\bar\caO_\bt$ we can give an upper bound to the number
of isolated branch points it can possess, using   
the intersection number $\Sigma_0\cdot\Sigma_\infty$. Since $\Gamma^*_{\beta_0}$ is invariant under 
$\gamma\mapsto -\gamma$ the
set of indices $\{l_j,m_j:j=1,\ldots,N\}$ in \eqref{eq:ellc} is also invariant under change of sign. 
The degree of the map $\ell^\C$ restricted to the first factor $\CP^1\times\{\text{pt}\}$ is,
generically, $\max_{j,m}|l_j-l_m|$. Since the $\gamma_j$ lie on a circle this equals $2l$, where
$l=\max_j\{l_j\}$. Similarly $\ell^\C$ has degree $2m$ where $m=\max_j\{m_j\}$ when restricted to the
second factor. Hence $\Sigma_0$ and $\Sigma_\infty$ are curves of type $(2l,2m)$ on $\CP^1\times\CP^1$
and have intersection number $8lm$. 
\end{rem}

\subsection{The moduli space of HSL tori.}
Here we will give a simple description of the moduli space $\caM(\Gamma,\beta_0)$
of all weakly conformal \HSL tori $f:\C/\Gamma\to\R^4$, based at $f(0)=0$ and with Maslov class
$\beta_0$, modulo base point preserving symplectic isometries and dilations of $\R^4$. 
As in the previous section, we still assume 
$N=\#\Gamma^*_{\beta_0}$ for each $j=1,\ldots,N/2$.
We have seen that $x(\lambda,0)$ is uniquely determined by its values at the roots of
$p(\lambda)$, and by the symmetries \eqref{eq:xsymmetry} it suffices to specify the $N/2$ 
vectors
\begin{equation}\label{eq:w}
w_j = \frac{1}{\sqrt{p(0)s_j^{N+1}}}x(s_j,0)\in\R^4,\ j=1,\ldots,N/2.
\end{equation}
It is convenient now to identify $\R^4$ with $\H$ so that $J$ represents left multiplication by
$i\in\H$ and $L$ represents left multiplication by $j\in\H$. In that case the action of $G_0$
corresponds to the action by right multiplication of the group $\Spin(3)$ of unit quaternions. It
follows from lemma \ref{lem:rotations} that the assignment
\begin{equation}\label{eq:wtof}
(w_1,\ldots,w_{N/2})\mapsto x(\lambda,0)\mapsto f
\end{equation}
given by \eqref{eq:fevolution} intertwines the action of $\H^*$ (by right multiplication) with the action
of $\R^+\times G_0$ (by dilations and symplectic isometries). 

We may allow $w_j$ to take any value in $\H$ provided they are not all zero: the subset on which they
are non-zero will correspond to $\Delta_f$. So \eqref{eq:wtof} descends to a bijection. 
In fact if we give the latter the manifold structure it inherits from the Banach
space $C^2(S^1\times S^1,\R^4)$, equipped with norm of uniform convergence on derivatives up to second
order, we see this must be a diffeomorphism. 
\begin{prop}\label{prop:moduli} 
Set $n=\frac{N}{2}-1$. Whenever $N>0$, the map 
$\HP^n\to \caM(\Gamma,\beta_0)$ defined by \eqref{eq:wtof} is a diffeomorphism.
\end{prop}
Of course, the two moduli spaces $S(\Gamma,\beta_0)$ and $\caM(\Gamma,\beta_0)$ are not isomorphic,
since the first labels $T_0$-orbits and the second labels $G_0$-orbits. Hence the former will be a
$G_0/T_0\simeq S^2$ bundle over the latter. The next result explains exactly what this is.
\begin{prop}
The action of $G_0$ on $\caS(\Gamma,\beta_0)$ via base point preserving isometries 
is equivalent, in appropriate homogeneous coordinates, to the action of $SU(2)$ on 
$\CP^{N-1}$ whose orbits are the fibres of $\CP^{2n+1}\to\HP^n$. 
\end{prop}
\begin{proof}
If we take the coordinates on $\caM(\Gamma,\beta_0)$ from $\eqref{eq:wtof}$ then we must make a 
change coordinates on the real manifold $\caS(\Gamma,\beta_0)$: we identify
the spectral data for $f$ with the point
\begin{equation}\label{eq:Lcoords}
[\chi_1(s_1,0),\chi_2(s_1,0),\ldots,\chi_1(s_{n+1},0),\chi_2(s_{n+1},0)]\in\CP^{2n+1}
\end{equation}
Here we are using the fact, from \eqref{eq:chi}, that 
$\chi_2(s_j,0) = is_j^{N+1}p(0)\overline{\chi_1(s_{j+N/2},0)}$ for
$1\leq j\leq n+1$. On the other hand the point in $\caM(\Gamma,\beta_0)$ corresponding to $f$ is
identified by the quaternionic homogeneous coordinates
\[
[w_1,\ldots,w_{n+1}]_\H\in\HP^n,\quad w_j = \frac{1}{\sqrt{p(0)s_j^{N+1}}}x(s_j,0).
\]
Here the identification of $\R^4$ with $\H$ is via
\[
w_j = \sum_{k=1}^4w_{jk}\e_k = (w_{j1} + w_{j2}J+w_{j3}L+w_{j4}JL)\e_1.
\]
Now we observe that 
\[
w_j =  \frac{1}{\sqrt{2}}( (w_{j1}+iw_{j2})v_1 + (w_{j3}+w_{j4})v_2 + 
(w_{j1}-iw_{j2})v_3 + (w_{j3}-w_{j4})v_4 ),
\]
and therefore
\[
\chi_1(s_j,0) = \sqrt{\frac{s_j^{N+1}p(0)}{2}}(w_{j1}+iw_{j2}),\quad
\chi_2(s_j,0) = \sqrt{\frac{s_j^{N+1}p(0)}{2}}(w_{j3}+iw_{j4}).
\]
Therefore the natural fibration $\CP^{2n+1}\to\HP^n$ obtained by the identification $\H=\C+\C j$ matches
the map $\caS(\Gamma,\beta_0)\to\caM(\Gamma,\beta_0)$ in the chosen coordinates.
\end{proof}

\section{Higher flows and Hamiltonian variations.} 

The action of the real subgroup $J_R\subset \Jac(X)$ generates the so-called higher flows. 
These give Lagrangian variations for our immersed Lagrangian tori, and 
we will show that a codimension $1$ subspace of these are actually Hamiltonian variations.
First we need to explain what we mean by Lagrangian and Hamiltonian vector fields along a 
Lagrangian immersion.

For a Lagrangian immersion
$f:M\to N$ into a \Kah manifold $N$ a section $T\in\Gamma(f^{-1}TN)$ is called a 
Lagrangian vector field \emph{along} $f$ if $\sigma_T=
f^*(T\rfloor\omega)$ is a closed 1-form, and $T$ is Hamiltonian when $\sigma_T$ is exact. We denote 
the vector space of Lagrangian vector fields along $f$ by $\caX_{\Lag}(f)$ and use $\caX_{\Ham}(f)$ 
to denote the subspace of Hamiltonian vector fields.
As one expects, families of Lagrangian immersions give rise to Lagrangian vector fields along an immersion \cite{CheM}. 

Now we consider the Lagrangian variations corresponding to the higher flows. For simplicity set
$\caS=\caS(\Gamma,\beta_0)$. 
Let $\caH\simeq\C^N\setminus\{0\}$ be the space of HSL tori in $\R^4$ 
given by \eqref{eq:fconstruction}, \eqref{eq:chi1} and \eqref{eq:chi-1} parameterised
by $\bt\in\C^N\setminus\{0\}$. The map $\C^N\to\CP^{N-1}$ which sends $\bt$ to $[\bt]$ is then
the parameterisation of the map $\caH\to\caS$ which assigns to each HSL torus its spectral data. 
Suppose $f\in\caH$ corresponds to $\bt= (t_1,\ldots,t_{N})$.
To each $\ba=(a_j)\in\C^{N}$ we can assign the Lagrangian vector field 
$T=(\partial f/\partial t)_{t=0}$ along $f$ corresponding to the Lagrangian deformation 
$f(t)$ determined by the real curve
\begin{equation}\label{eq:T}
t\mapsto (t_1e^{ta_1},\ldots,t_{N}e^{ta_{N}}),\quad t\in\R,
\end{equation}
in $\caS$. This gives a linear map
\begin{equation}\label{eq:Lag}
\caT_f:\C^N\to\caX_{\Lag}(f);\quad \caT_f(\ba) = T.
\end{equation}
It induces a linear map $\C^N/\C\to T_{[\bt]}\caS$ which is tangent to the action of $(\Ct)^N/\Ct$
which generates the higher flows for $f$ \eqref{eq:higherflows}. Notice that what we lose in the
quotient are all dilations and the rotations from $T_0$.

Now recall that $\caZ^1(\C/\Gamma)$ carries the Hodge inner product
\[
(\sigma_1,\sigma_2) = \frac{1}{A(\C/\Gamma)}\int_{\C/\Gamma}\sigma_1\wedge *\sigma_2,
\]
from which we obtain the orthogonal decomposition $\caZ^1(\C/\Gamma) = 
\caH^1(\C/\Gamma)\oplus\caB^1(\C/\Gamma)$, where the first summand is the space of harmonic 
1-forms. Since $d\beta$ and $*d\beta$ span $\caH^1(\C/\Gamma)$ we know that 
$\sigma_T\in\caZ^1(\C/\Gamma)$ is exact, and thus $T$ is Hamiltonian, precisely when
\[
(\sigma_T,d\beta)=0 = (\sigma_T,*d\beta).
\]
Notice that we can write
\[
\sigma_T = JT\cdot df,
\]
using dot product notation for the metric on $\R^4$ and its complex bilinear extension.
\begin{thm}\label{thm:higher}
Let $T=\caT_f(\ba)\in\caX_{\Lag}(f)$ be given by \eqref{eq:Lag}. Then 
$T$ is Hamiltonian along $f$ precisely when both real equations 
\begin{equation}\label{eq:c_0=0}
\sum_{j=1}^{N}\frac{\Re(a_j)|t_j|^2}{\prod_{k\neq j } |s_k-s_j|^2 } =0,\qquad
\sum_{j=1}^{N}\frac{\Re(s_j)}{\Im(s_j)}\frac{\Re(a_j)|t_j|^2}{\prod_{k\neq j } |s_k-s_j|^2 } =0.
\end{equation}
are satisfied.
\end{thm}
The first equation in \eqref{eq:c_0=0} is the statement that $\caT_f(\ba)$ is a stationary variation for
area (cf.\ the area formula \eqref{eq:area} below). Together, equations \eqref{eq:c_0=0} can also be
derived as the
condition that the Liouville form, which when pulled back along $f$ is the $1$-form 
$\sigma_f=Jf\cdot df$ on $\C/\Gamma$, has stationary $\Gamma$-periods for the variation $\caT_f(\ba)$.
These periods are well-known to be Hamiltonian isotopy invariants.
\begin{proof}
Let us begin by writing $\sigma_T=cdz+\bar c d\bar z$, $c= JT\cdot f_z$.
Since we are working over a torus, we can simplify computations significantly by observing 
that the orthogonal projection of $\sigma_T$ onto $\caH^1(\C/\Gamma)$ is just its Fourier zero 
mode, i.e., the coefficient of $e_0$ in the Fourier decomposition. 
Hence, if $c_0$ denotes the 
Fourier zero mode of $c$ then
\begin{eqnarray}\label{eq:hodge}
(\sigma_T,d\beta)& = &\pi(c_0dz+\bar c_0d\bar z,\bar\beta_0dz+\beta_0d\bar z)=
4\pi\Re(c_0\beta_0),  \\
(\sigma_T,*d\beta)&=&\pi(c_0dz+\bar c_0d\bar z,-i\bar\beta_0dz+i\beta_0d\bar z)=
-4\pi\Im(c_0\beta_0).
\end{eqnarray}
To compute $c_0$ we first simplify the formula \eqref{eq:fconstruction} by writing it as
\[
f = Av_1 +\bar A \bar v_1 + \bar Bv_2 + B \bar v_2,
\]
using
\[
A = \frac{1}{p(1)}(\chi_1(1,z)e_{\beta_0/2} -\chi_1(1,0)),\quad
B = \frac{1}{ip(1)}(\chi_1(-1,z)e_{-\beta_0/2} -\chi_1(-1,0)).
\]
Since $f(0)=0$ we may write each of the functions $A,B:\C/\Gamma\to\C$ as
\[
A = \sum_{j=1}^NA_j(e_{-\gamma_j+\beta_0/2}-1),
\quad B = \sum_{j=1}^NB_j(e_{-\gamma_j-\beta_0/2}-1),
\]
and, using equations \eqref{eq:chi1} and \eqref{eq:chi-1}, 
\[
A_j = \frac{(-1)^{j-1} }{ p(1)} \Delta_j t_j , \quad 
B_j = \frac{(-1)^{j-1} }{ i p(1)} \tilde{\Delta}_j t_j.
\] 
The deformation $f(t)$ is given by
\[
A(t) = \sum_{j=1}^N e^{ta_j}A_j(e_{-\gamma_j+\beta_0/2}-1),\quad
B(t) = \sum_{j=1}^N e^{ta_j}B_j (e_{-\gamma_j-\beta_0/2}-1).
\]
Since $v_1,\bar v_1,v_2,\bar v_2$ is a Hermitian orthonormal
frame for $\C^4$ we find that 
\[
Jf_t\cdot f_z = i[\frac{\partial A}{\partial t}\frac{\partial\bar A}{\partial z} - 
\frac{\partial\bar A}{\partial t}\frac{\partial A}{\partial z} - 
\frac{\partial B}{\partial t}\frac{\partial\bar B}{\partial z} + 
\frac{\partial\bar B}{\partial t}\frac{\partial B}{\partial z} ].
\] 
We compute at $t=0$ 
\begin{eqnarray*}
\partial A/\partial t=\sum_{j=1}^{N}  a_j A_j (e_{\beta_0/2-\gamma_j}-1),
&\partial \bar{A}/\partial t=\sum_{j=1}^{N}  \bar{a}_j \bar{A}_j (e_{\gamma_j-\beta_0/2}-1)\\
\partial A/\partial z = i \pi \sum_{j=1}^{N}  \left(\frac{\bar{\beta}_0}{2}-\bar{\gamma}_j \right) 
A_j e_{\beta_0/2-\gamma_j},
&\partial \bar{A}/\partial z = -i\pi\sum_{j=1}^{N}\left(\frac{\bar{\beta}_0}{2}-\bar{\gamma}_j \right) 
\bar A_j e_{\gamma_j-\beta_0/2}
\end{eqnarray*}
and similar expressions for $B$. The Fourier zero mode $c_0$ of $JT\cdot f_z$ is given by
the following calculation. To simplify it, we let $V_N$ stand for the Vandermonde determinant
$\prod_{1\leq k<\ell\leq N}(s_\ell-s_k)$ (cf.\ appendix \ref{app:theta} below).
\begin{eqnarray*}
c_0 &=& \pi \bigg( 
\sum_{j=1}^{N}   \left(\frac{\bar{\beta}_0}{2}-\bar{\gamma}_j \right) (a_j+\bar{a}_j) |A_j|^2
-
\sum_{j=1}^{N}   \left(-\frac{\bar{\beta}_0}{2}-\bar{\gamma}_j \right) (a_j+\bar{a}_j) |B_j|^2
\bigg)\\
&=& 2 \pi \sum_{j=1}^{N}  \bigg( 
\left( \frac{\bar{\beta}_0}{2}-\bar{\gamma}_j \right) |A_j|^2
+ \left(\frac{\bar{\beta}_0}{2}+\bar{\gamma}_j \right) |B_j|^2
\bigg) \Re (a_j) \\
&=& \frac{ \pi \bar{\beta}_0 }{ |p(1)|^2 }  \sum_{j=1}^{N} \bigg( 
\left( 1-\bar{s}_j \right) |\Delta_j|^2 + \left(1+\bar{s}_j \right) |\tilde{\Delta}_j|^2
\bigg) \Re (a_j)|t_j|^2 \\
&=&  \pi \bar{\beta}_0 |V_N|^2 \sum_{j=1}^{N} \left( \frac{1}{1-s_j} + \frac{1}{1+s_j} \right)
\frac{ \Re (a_j)|t_j|^2 }{ \prod_{k\neq j } |s_k-s_j|^2 } 
\end{eqnarray*}
where we have used the property $p(1)=p(-1)$ to write 
\[
|\Delta_j|^2=\frac{ |p(1)|^2 |V_N|^2}{|1-s_j|^2 \prod_{k\neq j } |s_k-s_j|^2 }
= \frac{ |1+s_j|^2 }{ |1-s_j|^2 } |\tilde{\Delta}_j|^2 .
\]
It follows that
\begin{eqnarray}
(\sigma_T,d\beta) & = & 4\pi^2|\beta_0|^2|V_N|^2 \sum_{j=1}^{N}\frac{ \Re (a_j)|t_j|^2 }{
\prod_{k\neq j } |s_k-s_j|^2 }\\
(\sigma_T,*d\beta) & = & -4\pi^2|\beta_0|^2|V_N|^2 \sum_{j=1}^{N}\frac{\Re(s_j)\Re (a_j)|t_j|^2 }{
\Im(s_j)\prod_{k\neq j } |s_k-s_j|^2 }.
\end{eqnarray}
\end{proof}
\begin{rem}\label{rem:area}
It is interesting to note that, given the isomorphism $J_R\simeq (S^1)^{N-1}\times(\R^+)^{N-1}$,
all the higher flows tangent to the compact factor of $J_R$ are Hamiltonian, but not \emph{all}
higher flows are Hamiltonian, in the sense that there are vector fields $\caT_f(\ba)$ which are not
Hamiltonian and whose projections onto $\caS$ give non-trivial higher flows . This is a reflection 
of the fact that the area functional is non-constant on $\caS$.
Using the expressions above it is straightforward to show that the area of $f(\C/\Gamma)$, when $f$
corresponds to $[\bt]\in\CP^{N-1}$, is given by
\begin{equation}\label{eq:area}
A(f)=\int_{\C/\Gamma} |f_z|^2|dz|^2 = A(\C/\Gamma) \frac{\pi^2|\beta_0|^2|V_N|^2}{2}
\sum_{j=1}^{N} \frac{|t_j|^2}{\prod_{k\neq j } |s_k-s_j|^2 }.
\end{equation}
Here, of course, we run into the
problem that the dilations act on the area but have been factored out of the space 
$\caS$ but we can restrict $A(f)$ to the unit sphere
$S^{2N-1}\subset\C^N$: it then descends to $\CP^{N-1}$. It is tantalising to observe that,
up to a constant independent of $\bt$, $A(f)$ equals
\begin{equation}
\int_{\C/\Gamma}|\theta_\infty(\bt(z))|^2|dz|^2,
\end{equation}
where $\theta_\infty$ is defined in \eqref{eq:thetai}.
\end{rem}
\begin{rem}\label{rem:genus}
Our final remark concerns the dependence of $A(f)$ on the spectral genus of $f$. 
Here we interpret the spectral genus of $f$ to be the arithmetic genus $g_f$ of $X_f$. From the earlier
discussion $g_f=4(N_1-1)$ where $N_1$ is the complex dimension of the $J_R$-orbit of the spectral data
for $f$ in $\caS(\Gamma,\beta_0^*)$. The closure of these $J_R$-orbits provides a stratification of 
$\caS(\Gamma,\beta_0^*)$, and as we have observed this agrees with the stratification of $\CP^{N-1}$ 
by $(\Ct)^N/\Ct$-orbits. In particular, one can decrease the spectral genus by going to the boundary of
a $J_R$-orbit. In our paramaterisation this corresponds to taking the limit as some $t_j\to 0$.

To make proper sense of the dependence of $A(f)$ on spectral genus we can restrict $A(f)$ to the unit 
sphere $S^{2N-1}\subset\C^N$: it then descends to $\CP^{N-1}$. Clearly ``going to the boundary of
spectral genus'' decreases $A(f)$. Moreover,
it is easy to show, using Lagrange multipliers, that $A(f)$ has critical points
only at the homogeneous tori (those which arise as orbits of a homomorphism from $\C/\Gamma$ into
$G$). These are all minima of $A(f)$. These are also the HSL tori with 
lowest (non-trivial) spectral genus. While both of these observations are elementary for this geometry,
they are an echo of what one hopes to find in other, more complicated, integrable surface theory (cf.
the work of Kilian and Schmidt \cite{KilS} on deformations of CMC tori into cylinders, which is a form
of ``going to the boundary of spectral data'' to decrease spectral genus, for example).
\end{rem}

\appendix

\section{The $\theta$-function.}\label{app:theta}

Here we will prove proposition \ref{prop:varphi} and also explain its geometric meaning: it is
essentially the Riemann vanishing theorem for a singular curve related to our spectral curve $X$.

At an algebraic level the proof is straightforward. For any $h\in F_{N-1}$ write
$h(\lambda)=\sum_{k=0}^{N-1}h_k\lambda^k$.
Given $[\bt]=[t_1,\ldots,t_N]\in\CP^{N-1}$ we have $\iota([h])=[\bt]$ precisely when we can solve the
linear system 
\begin{equation}\label{eq:vandermonde}
\begin{pmatrix}
1 & s_1 & \ldots & s_1^{N-1}\\
\vdots & \vdots & &\vdots\\ 
1& s_{N-1} &\ldots & s_{N-1}^{N-1}\\
1 & s_N&\ldots&s_N^{N-1}
\end{pmatrix}
\begin{pmatrix} h_0\\ \vdots \\h_{N-2} \\ h_{N-1}\end{pmatrix}
= \begin{pmatrix} t_1\\ \vdots \\ t_{N-1} \\ t_N\end{pmatrix}
\end{equation}
The matrix has determinant $V_N=\prod_{m>j}(s_m-s_j)$ (it is a Vandermonde matrix), hence it is
invertible whenever $s_1,\ldots,s_N$ are distinct. 
Hence $\iota$ is an isomorphism. Further, $\iota([\varphi(\lambda,\bt)])=[\bt]$ since
\[
\varphi(s_j,\bt) = 
(\prod_{m\neq j}(s_m-s_j))\det\begin{pmatrix}
1 & s_1 & \ldots & s_1^{N-2} & 0\\
\vdots & \vdots & &\vdots& \vdots\\ 
1 & s_j&\ldots&s_j^{N-2} & t_j\\
\vdots & \vdots & &\vdots& \vdots\\ 
1 & s_N&\ldots&s_N^{N-2} & 0
\end{pmatrix}
=
t_j V_N.
\]
Finally, by applying Cramer's rule to
\eqref{eq:vandermonde} it follows that $h_{N-1}=0$ precisely when 
\begin{equation}\label{eq:thetai}
\theta_\infty(\bt)=
\det\begin{pmatrix}
1 & s_1 & \ldots & s_1^{N-2} & t_1\\
\vdots & \vdots & &\vdots& \vdots\\ 
1 & s_N&\ldots&s_N^{N-2} & t_N
\end{pmatrix}= V_N\sum_{j=1}^N(\prod_{m\neq j}\frac{1}{s_m-s_j})t_j=0.
\end{equation}
Similarly, $h_0=0$ when
\begin{equation}\label{eq:theta0}
\theta_0(\bt)=
\begin{pmatrix}
t_1 & s_1 & \ldots & s_1^{N-1}\\
\vdots & \vdots & &\vdots\\ 
t_N & s_N&\ldots&s_N^{N-1}
\end{pmatrix}
=(\prod_{j=1}^Ns_j)V_N\sum_{j=1}^N(\prod_{m\neq j}\frac{1}{s_m-s_j})\frac{t_j}{s_j}=0.
\end{equation}
This completes the proof.

The geometric meaning is as follows (cf.\ the discussion of rational nodal curves in \cite[3.251]{Mum}). 
Let $Y$ be the rational singular curve obtained from $\hat\C$ by
identifying the points $s_1,\ldots,s_N$ into one singularity $s$. This curve has arithmetic genus
$g=N-1$ and one can show that $\Jac(Y)\simeq
(\Ct)^{N-1}$. It has an Abel map
\[
\caA_\infty:Y\setminus\{s\}\to\Jac(Y);\quad Q\mapsto \caO_Y(Q-\infty).
\]
We can identify $Y\setminus\{s\}$ with $\hat\C\setminus\fS$ and $\Jac(Y)$ with 
\[
\{[t_1,\ldots,t_{N-1},t_N]\in\CP^{N-1}:t_j\neq 0 \forall j\}.
\]
Since $\caO_Y(Q-\infty)$ has unique, up to scale, global section $\lambda-Q$ the Abel map
extends holomorphically to the map 
\[
\caA_\infty:\hat\C\to\CP^{N-1};\quad \caA_\infty(Q) =
[s_1-Q,\ldots,s_N-Q].
\]
By analogy with the case of compact Riemann surfaces, we say the $\theta$-divisor for $Y$ (given base
point $\infty$) is the image of $(Y\setminus\{s\})^{(N-2)}$ in $\Jac(Y)$ under the Abel map
on divisors (note that $N-2=g-1$). 
This image is the set of all divisor classes of the form $[E-(N-2).\infty]$, hence corresponds to
polynomials of degree strictly less than $N-1$. Therefore this $\theta$-divisor is 
$\Jac(Y)\cap\Theta_\infty$. Similarly, $\Theta_0$ is a translate of this $\theta$-divisor.
The function $\varphi$ plays the role of a translate of the $\theta$-function pulled back to 
$Y$ along the Abel map. In this context proposition \ref{prop:varphi} is the analogue of Riemann's 
vanishing theorem.


The link with our spectral curve $X$ is that there is a isomorphism of real groups $\Jac(Y)\simeq J_R$
given by identifying $Y\setminus\{s\}$ with $C_1\setminus\fS$ and mapping
\[
\caO_Y(D)\to\caO_X(D+\rho^*D+\tau^*D+(\rho\tau)^*D),\quad D\in\Div_0(Y\setminus\{s\}).
\]


\begin{thebibliography}{444}

\bibitem{Anc} H Anciaux, {\em An isoperimetric inequality for Hamiltonian stationary Lagrangian tori in
$\C^2$ related to Oh's conjecture.} Math.\ Z.\ 241 (2002), no.\ 3, 639--664.

\bibitem{Anc1} H Anciaux, {\em Construction of many Hamiltonian stationary Lagrangian surfaces in
Euclidean
four-space.}  Calc.\ Var.\ Partial Differential Equations  17  (2003),  no.\ 2, 105--120.

\bibitem{Anc2} H Anciaux, {\em Hamiltonian stationary Lagrangian surfaces in Euclidean four-space
andrelated minimization problems. II. Some minimization results.}  Ann.\ Global Anal.\ Geom.\  22
(2002), no.\ 4, 341--353.

\bibitem{Boh08} C Bohle, \textit{Constrained Willmore tori in the $4$-sphere.}, arXiv 0803.06331v1.

\bibitem{BohLPP} C Bohle, K Leschke, F Pedit \& U Pinkall, \textit{Conformal
maps from a 2-torus to the $4$-sphere}, arXiv 0712.2311v1 (2007).

\bibitem{BurK} F E Burstall \& I Khemar, \textit{Twistors, $4$-symmetric spaces and integrable
systems}, Math.\ Ann.\  344  (2009), 451--461.

\bibitem{CheM} B-Y Chen \& J-M Morvan, \textit{Deformations of isotropic submanifolds in K\" ahler manifolds}. J.\ Geom.\ Phys.\ 13 (1994), 79--104.

\bibitem{HelR} F H\' elein \& P Romon, \textit{Hamiltonian stationary Lagrangian
surfaces in $\C^2$}, Comm.\ Anal.\ Geom.\ 10 (2002), 79--126.

\bibitem{HelR2} F H\' elein \& P Romon, {\em Hamiltonian stationary Lagrangian surfaces in Hermitian
symmetric spaces.}  Differential geometry and integrable systems (Tokyo, 2000),  161--178, Contemp.\
Math., 308, Amer.\ Math.\ Soc., Providence, RI, 2002.

\bibitem{Ilm}  T Ilmanen, {\em A note of the Hamiltonian area conjecture.} ETH Zurich preprint 1998. 

\bibitem{KilS} M Kilian \& M U Schmidt, \textit{On the moduli of constant mean curvature cylinders
of finite type in the $3$-sphere}, arXiv 0712.0108v2 (2008).

\bibitem{LesR} K Leschke \& P Romon, \textit{Darboux transforms and spectral curves of
Hamiltonian stationary Lagrangian tori.} Calc.\ Var.\ and PDE (2009) DOI 10.1007/s00526-009-0278-6. 

\bibitem{McI} I McIntosh, \textit{Harmonic tori and their spectral data.}
``Surveys on Geometry and Integrable Systems'', Adv.\ Studies Pure Math.\ 51 (2008).

\bibitem{Mor} J-M Morvan, \textit{Classe de Maslov d'une immersion lagrangienne et minimalit\' e.} 
C.\ R.\ Acad.\ Sci.\ Paris Sér.\ I Math.\  292  (1981), no. 13, 633--636.

\bibitem{Mum} D Mumford, \textit{Tata lectures on Theta II; Jacobian theta functions and differential
equations.} With the collaboration of C Musili, M Nori, E Previato, M Stillman and H Umemura.
Progress in Mathematics, 43.\ Birkh\" auser Boston, Inc., Boston, MA, 1984.

\bibitem{Oh2} Y G Oh, {\em Volume minimization of Lagrangian submanifolds under Hamiltonian
deformations.}
Math.\ Z.\ 212, 175--192 (1993).

\bibitem{Ser} J-P Serre, {\em Algebraic groups and class fields,} Graduate
Texts in Math.\ 117, Springer, New York (1988).

\bibitem{SchW} R Schoen \& J Wolfson, {\em Minimizing volume among Lagrangian submanifolds.} Proc.\
Symp.\ Pure Math.\ 65, 181--199 (1999).

\bibitem{SchW1} R Schoen \& J Wolfson, {\em Minimizing area among Lagrangian surfaces: the mapping
problem.} J.\ Differential Geom.\ 58 (2001), no.\ 1, 1--86.


\end{thebibliography}
\end{document}